\documentclass[12pt]{article}

\usepackage{amssymb,amsmath,exscale}

\usepackage{graphicx}

\usepackage{color}
\definecolor{marin}{rgb}   {0.,   0.3,   0.7} 
\definecolor{rouge}{rgb}   {0.8,   0.,   0.} 
\definecolor{sepia}{rgb}   {0.8,   0.5,   0.} 
\usepackage[colorlinks,citecolor=marin,linkcolor=rouge,
            bookmarksopen,
            bookmarksnumbered
           ]{hyperref}

\newtheorem{lemma}{Lemma}[section]

\newtheorem{theorem}[lemma]{Theorem}
\newtheorem{proposition}[lemma]{Proposition}

\newtheorem{remark}[lemma]{Remark}
\newtheorem{example}[lemma]{Example}
\newtheorem{hypothesis}[lemma]{Hypothesis}
\newtheorem{notation}[lemma]{Notation}
\newtheorem{definition}[lemma]{Definition}
\newtheorem{conclusion}[lemma]{Conclusion}
\numberwithin{equation}{section}
\newcommand{\QED}{\mbox{}\hfill \raisebox{-0.2pt}{\rule{5.6pt}{6pt}\rule{0pt}{0pt}} 
          \medskip\par}             
\newenvironment{Proof}{\noindent
    \parindent=0pt\abovedisplayskip = 0.5\abovedisplayskip
    \belowdisplayskip=\abovedisplayskip{\bfseries Proof. }}{\QED}
\newenvironment{Proofof}[1]{\noindent
    \parindent=0pt\abovedisplayskip = 0.5\abovedisplayskip
    \belowdisplayskip=\abovedisplayskip{\bfseries Proof of #1. }}{\QED}

\newcommand{\Bc}{\mathsf{N}}
\newcommand{\Rc}{\mathsf{R}}

\newcommand{\Ac}{\mathcal{A}}
\newcommand{\ab}{\boldsymbol{a}}
\newcommand{\bb}{\boldsymbol{b}}

\newcommand{\dd}{\mathrm{d}}

\newcommand{\Hc}{\mathcal{H}}

\newcommand{\Ic}{\mathcal{I}}
\newcommand{\Jc}{\mathcal{J}}

\newcommand{\jb}{{\boldsymbol{j}}}
\newcommand{\kb}{{\boldsymbol{k}}}

\newcommand{\N}{\mathbb{N}}
\newcommand{\Nc}{\mathcal{N}}

\newcommand{\Mc}{\mathcal{M}}
\newcommand{\Pc}{\mathcal{P}}

\newcommand{\R}{\mathbb{R}}

\newcommand{\C}{\mathbb{C}}

\newcommand{\T}{\mathbb{T}}
\newcommand{\Z}{\mathbb{Z}}
\newcommand{\Tc}{\mathcal{T}}

\newcommand{\Uc}{\mathcal{U}}

\newcommand{\Wc}{\mathcal{W}}

\newcommand{\Zc}{\mathcal{Z}}

\newcommand{\Norm}[2]{\|#1\|\left.\vphantom{T_{j_0}^0}\!\!\right._{#2}}         
\newcommand{\SNorm}[2]{|#1|\left.\vphantom{T_{j_0}^0}\!\!\right._{#2}}             

\title{Birkhoff normal form and splitting methods\\ for semi linear Hamiltonian PDEs. \\
Part II: Abstract splitting. }        
\author{Erwan Faou, Beno\^it Gr\'ebert and Eric Paturel}       
\begin{document}
\maketitle
\abstract{
We consider Hamiltonian PDEs that can be split into a linear unbounded operator and a regular non linear part. We consider abstract splitting methods associated with this decomposition where no discretization in space is made. We prove a normal form result for the corresponding discrete flow under generic non resonance conditions on the frequencies of the linear operator and on the step size. This result implies the conservation of the regularity of the numerical solution associated with the splitting method over arbitrary long time, provided the initial data is small enough. This result holds for numerical schemes controlling the round-off error at each step to avoid possible high frequency energy drift.  We apply this results to nonlinear Schr\"odinger equations as well as the nonlinear wave equation. }


\tableofcontents

\section{Introduction}

In this work, we consider a class of Hamiltonian partial differential equations whose Hamiltonian functions $H=H_0+P$ can be divided into a linear unbounded operator $H_0$ with discrete spectrum and a non linear function $P$ having a zero of order at least $3$ at the origin of the phase space. Typical examples are given by the non linear wave equation  or the non linear Schr\"odinger equation on the torus. 

Amongst all the numerical schemes that can be applied to these Hamiltonian PDEs, splitting methods entail many advantages, as they provide symplectic and explicit schemes, and can be easily implemented using fast Fourier transform if the spectrum of $H_0$ expresses easily in  Fourier basis. Generally speaking, a splitting schemes is based on the approximation 
\begin{equation}
\label{E001}
\varphi_H^h \simeq \varphi_{H_0}^h \circ \varphi_{P}^h
\end{equation}
for small time $h$, and where $\varphi_{K}^t$ denotes the exact flow of the Hamiltonian PDE associated with the Hamiltonian $K$. 

The understanding of the long-time behavior of splitting methods for Hamiltonian PDEs is a fundamental ongoing challenge in the field of geometric integration, as the classical arguments of {\em backward error analysis} (see for instance \cite{HLW}) do not applied in this situation where the frequencies of the system are arbitrary large, and where resonances phenomenon are known to occur for some values of the step size. Recently, many progresses have been made in this direction. A first result using normal form techniques was given by {\sc Dujardin \& Faou} in \cite{DF07} for the case of the linear Schr\"odinger equation with small potential. Concerning the non linear case, results exists by {\sc Cohen, Hairer \& Lubich}, see \cite{CHL08b,CHL08c}, for the wave equation and {\sc Gauckler \& Lubich}, see \cite{GL08a,GL08b}, for the nonlinear Schr\"odinger equation using the technique of modulated Fourier expansion. However to be valid these results use non-resonance conditions that are generically satisfied only under CFL conditions linking the step-size $h$ and highest frequencies of the space discretization of the Hamiltonian PDE. 

In this paper, we use a normal form techniques to prove the long-time preservation of the regularity of the initial solution under {\em generic} non resonances conditions valid for a large set of equations and a large set of time step, independently of the space discretization parameter.  

Normal form techniques have proven to be one of the most important tool for the understanding of the long time behaviour of Hamiltonian PDE (see \cite{Bam03, BG06, Greb07, Bam07, BDGS, GIP}). Roughly speaking, the dynamical consequences of such results are the following: starting with a small initial value of size $\varepsilon$ in a Sobolev space $H^s$, then the solution remains small in the same norm over long time, namely for time $t \leq C_r \varepsilon^{-r}$ for arbitrary $r$ (with a constant $C_r$ depending on $r$). Such results hold under {\em generic} non resonance conditions on the frequencies of the underlying linear operator $H_0$ associated with the Hamiltonian PDE, that are valid in a wide number of  situations (nonlinear Schr\"odinger equation on a torus of dimension $d$ or with Dirichlet boundary conditions, nonlinear wave equation with periodic or Dirichlet conditions in dimension 1, Klein Gordon equation on spheres or Zoll manifolds.). 

In this paper, we mainly show that the same kind of results hold true for {\em numerical} solutions associated with the abstract splitting method \eqref{E001} under some further restrictions specifically induced by the time discretization.

This work is the second of a series of two.

In the first part \cite{FGP1}, we consider full discretizations of the Hamiltonian PDE, with a spectral discretization parameter $K$. We show that under the hypothesis $K \leq \varepsilon^{-\sigma}$ for some constant $\sigma$ depending on the precision degree $r$, the same conclusion as in the continuous case can be drawn. Though concerning the discretization of a Hamiltonian PDE, the method used in \cite{FGP1} is by essence a finite dimensional Birkhoff normal form result using techniques that are rather classic in the dynamical system world. 

In some sense, the present paper studies the case where $K > \varepsilon^{-\sigma}$ by considering the splitting method where no discretization in space is made (i.e. $K = +\infty$). The techniques used involve the abstract framework developed in \cite{BG06, Greb07, Bam07}. However, instead of being valid for the (exact) abstract splitting \eqref{E001}, we have to consider {\em rounded} splitting methods of the form
\begin{equation}
\label{E002}
\Pi_{\eta,s}\circ \varphi_{H_0}^h\circ \varphi_P^h
\end{equation}
where $\Pi_{\eta,s}$ puts to zero all the frequencies $\xi_j$ whose weighted energy $|j|^{2s}|\xi_j|^{2s}$ in the Sobolev space $H^s$ is smaller than a given threshold $\eta^2$. Hence, for small $\eta$, \eqref{E002} is very close to the exact splitting method \eqref{E001}. The good news is that this threshold can be taken of the order $\varepsilon^r$, making this projection $\Pi_{\eta,s}$ very close to the identity, and in any case producing an error that is far beyond the round-off error in a computer simulation (particularly for large $s$).  

Our main result is given by Theorem \ref{Tmain}. 

\section{Description of the method}

Before going on into the precise statements and proofs, we would like to give  tentative explanations of the modifications observed in comparison with the continuous case. 

The method used in \cite{BG06} to prove the long-time conservation of Sobolev norms for small data is to start from a Hamiltonian $H = H_0 + P$ depending on an infinite number of variable $(\xi_j,\eta_k)$, $j, k \in \N$, and for a fixed number $r$, to construct a hamiltonian transformation $\tau$ close to the identity, and such that in the new variable, the Hamiltonian can be written
\begin{equation}
\label{EnfBG}
H_0  + Z + R
\end{equation}
where $Z$ is a real Hamiltonian depending only on the action $I_j = \xi_j\eta_j$ and $R$ a real Hamiltonian having a zero of order $r$. 

The key for this construction is an induction process where, at each step, the solution of an homological equation of the form 
\begin{equation}
\label{HomBG}
\{ H_0, \chi\} + Z = G
\end{equation}
where $G$ is a given homogeneous polynomial of order $n$, and where $Z$ depending only on the action and $\chi$ are unknown. Assume that $G$ is of the form 
$$
G = G_{\jb\kb} \,\xi_{j_1}\cdots \xi_{j_p}\eta_{k_1}\cdots \eta_{k_q}
$$
where $G_{\jb\kb}$ is a coefficients, $\jb = (j_1,\ldots,j_p) \in \N^p$ and $\kb = (k_1,\ldots,k_q) \in \N^q$. 
Then it is easy to see that the equation \eqref{HomBG} can be written 
\begin{equation}
\label{HomBG2}
\Omega(\jb,\kb) \chi_{\jb\kb} + Z_{\jb\kb} = G_{\jb\kb}
\end{equation}
where 
$$
\Omega(\jb,\kb) = \omega_{j_1} + \cdots + \omega_{j_p} - \omega_{k_1} -\cdots - \omega_{j_q}
$$
and where $Z_{\jb\kb}$ and $G_{\jb\kb}$ are unknown coefficients. 

It is clear that for $\jb = \kb$ (up to a permutation), we have $\Omega(\jb,\kb) = 0$ which imposes $Z_{\jb\kb} = G_{\jb\kb}$. When $\jb \neq \kb$ (taking into account the permutation), the solution of \eqref{HomBG2} relies on a non resonance conditions on the small divisors $\Omega(\jb,\kb)^{-1}$. 

In \cite{BG06}, {\sc Bambusi \& Gr\'ebert} use a non resonance condition of the form 
\begin{equation}
\label{EnonresBG}
\forall\, \jb \neq \kb,\quad | \Omega(\jb,\kb) | \geq \gamma \mu(\jb,\kb)^{-\alpha}
\end{equation}
where $\mu(\jb,\kb)$ denotes the third largest integer amongst $|j_1|,\ldots,|k_q|$. They moreover show that such a condition is guaranteed in a large number of situations (see \cite{BG06}, \cite{Greb07} or \cite{Bam07} for precise results). 

 Considering now the splitting method $\varphi_{H_0}^h \circ \varphi_{P}^h$, we see that we cannot work directly at the level of the Hamiltonian. To avoid this difficulty, we embed the splitting into the family of applications
$$
[0,1] \ni \lambda\mapsto \varphi_{H_0}^h \circ \varphi_{hP}^\lambda
$$
and we derive this expression with respect to $\lambda$ to work in the tangent space, where it is much more easy to identify real Hamiltonian than unitary flows. 

This explains why we deal here with time-dependent Hamiltonian. Note that we do not expand the operator $\varphi_{H_0}^h$  in powers of $h$, as this would yields positive powers of the unbounded operator $H_0$ appearing in the series. Unless a CFL condition is employed, this methods do not give the desired results (and do not explain the resonance effects observed for some specific values of $h$, see \cite{FGP1}). 

Now, instead of \eqref{HomBG}, the Homological equation appearing for the splitting methods is given in a discrete form
\begin{equation}
\label{HomFGP}
\chi \circ \varphi_{H_0}^h - \chi + Z = G. 
\end{equation}
In terms of coefficients, this equations yields
$$
(e^{ih\Omega(\jb,\kb)} - 1) \chi_{\jb\kb} + Z_{\jb\kb} = G_{\jb\kb}.
$$
The main difference with \eqref{HomBG2} is that we have to avoid not only the indices $(\jb,\kb)$ so that $\Omega(\jb,\kb) = 0$, but all of those for which $h \Omega(\jb,\kb) = 2 m \pi$ for some (unbounded) integer $m$.

In the case of a fully discretized system for which $\nabla_{z_j} P \equiv 0$ for $|j| > K$, then under the CFL-like condition of the form $h K^2 \leq C$ for some constant $C$ depending on $r$, then we have $|h \Omega(\jb,\kb)| \leq  \pi$, and hence 
\begin{equation}
\label{EnonresGL}
|e^{ih\Omega(\jb,\kb)} - 1| \geq h \gamma \mu(\jb,\kb)^{-\alpha}
\end{equation}
\eqref{EnonresGL}
is then a consequence of \eqref{EnonresBG}. Under this assumption, we can apply the same techniques used in  \cite{BG06} and draw the same conclusions. This is the kind of assumption made in \cite{CHL08b} and \cite{GL08b}.

The problem with \eqref{EnonresGL} is that it is non generic in $h$ outside the CFL regime. For example, in the case of the Schr\"odinger equation, the frequencies of the operator $H_0$ are such that  $\omega_j \simeq j^2$. Hence, for large $N$, if $(j_1,\ldots,j_p,k_1,\ldots,k_q)$ is such that $j_1 = N +1$, $k_1 = N$ and all the other are of order $1$ ($N$ is large here), we have $\Omega(\jb,\kb) \simeq  (N+1)^2 - N^2 \simeq 2N$. Hence, 
$$
|e^{ih\Omega(\jb,\kb)} - 1| \simeq | e^{2ihN} - 1 | 
$$
cannot be assumed to be greater than $h \gamma \mu(\jb,\kb)^{-\alpha} \simeq h$. Note that a generic hypothesis on $h$ would be here that this small divisor is greater than $h \gamma N^{-\alpha}$ for some constants $\gamma$ and $\alpha$. This means that we cannot control the small divisors $|e^{ih\Omega(\jb,\kb)} - 1|$ associated with the splitting scheme by the {\em third largest} integer in the multi index, but by the {\em largest}. 

Using a generic condition on $h \leq h_0$, we can prove a normal form result and show that the flow is conjugated to the flow of a hamiltonian vector field of the form \eqref{EnfBG}, but where $Z$ now contains terms depending only on the actions, and supplementary terms containing  at least two large indices. Here, large means greater than $\varepsilon^{-\sigma}$ where $\sigma$ depends on $r$. 

In the case of a full discretization of the Hamiltonian PDE with a spectral discretization parameter $K$, we thus see that if $K \leq \varepsilon^{-\sigma}$ then the normal form term $Z$ actually depends only on the actions, as the high frequencies greater that $\varepsilon^{-\sigma}$ are not present. This is essentially the result of \cite{FGP1}. 

In the case where $K > \varepsilon^{-\sigma}$, the normal form result that we obtain can be interpreted as follows: the non conservation of the actions can only come from two high modes (of order greater than $\varepsilon^{-\sigma}$) interacting together and contaminating the whole spectrum. The role of the projection operator $\Pi_{\eta,s}$ is to destroy these high modes at each step but only when these high modes have an energy greater than $\eta$ (cf. \eqref{Pi-1}). As we can take $\eta = \varepsilon^r$, the error induced is very small, and in particular, far beyond the round-off error in the numerical simulation. Note that the complete understanding of the numerical phenomenon, and in particular the possible interaction (or not) between two high modes would in principle require the introduction of round-off effects by adding stochastic terms in high frequencies. This is clearly out of the scope of this paper.

\section{Setting of the problem}

\subsection{Abstract Hamiltonian formalism}

We denote $\Nc = \Z^d$ or $\N^d$ (depending on the concrete application) for some $d \geq 1$.  For $a = (a_1,\ldots,a_d) \in \Nc$, we set
$$
|a|^2 = \max\big(1,a_1^2 + \cdots + a_d^2\big). 
$$
We consider the set of variables $(\xi_a,\eta_b) \in \C^{\Nc} \times \C^{\Nc}$ equipped with the symplectic structure
\begin{equation}
\label{Esymp}
i \sum_{a \in \Nc} \dd \xi_a \wedge \dd \eta_a. 
\end{equation}
We define the set $\Zc = \Nc \times \{ \pm 1\}$. For $j = (a,\delta) \in \Zc$, we define $|j| = |a|$ and we denote by $\overline{j}$ the index $(a,-\delta)$. 

We will identify a couple $(\xi,\eta)\in \C^{\Nc} \times \C^{\Nc}$ with 
$(z_j)_{j \in \Zc} \in \C^{\Zc}$ via the formula
$$
j = (a,\delta) \in \Zc  \Longrightarrow 
\left\{
\begin{array}{rcll}
z_{j} &=& \xi_{a}& \mbox{if}\quad \delta = 1,\\[1ex]
z_j &=& \eta_a & \mbox{if}\quad \delta = - 1,
\end{array}
\right.
$$
By  a slight abuse of notation, we often write $z = (\xi,\eta)$ to denote such an element. 

For a given real number $s \geq 0$, we consider the Hilbert space $\Pc_s = \ell_s(\Zc,\C)$ made of elements $z \in \C^{\Zc}$ such that
$$
\Norm{z}{s}^2 := \sum_{j \in \Zc} |j|^{2s} |z_j|^2 < \infty,
$$
and equipped with the symplectic form \eqref{Esymp}. 

Let $\Uc$ be a an open set of $\Pc_s$. For a function $F$ of $\mathcal{C}^1(\Uc,\C)$, we define its gradient by 
$$
\nabla F(z) = \left( \frac{\partial F}{\partial z_j}\right)_{j \in \Zc}
$$
where by definition, we set for $j = (a,\delta) \in \Nc \times \{ \pm 1\}$, 
$$
 \frac{\partial F}{\partial z_j} =
  \left\{\begin{array}{rll}
 \displaystyle  \frac{\partial F}{\partial \xi_a} & \mbox{if}\quad\delta = 1,\\[2ex]
 \displaystyle \frac{\partial F}{\partial \eta_a} & \mbox{if}\quad\delta = - 1.
 \end{array}
 \right.
$$
Let $H(z)$ be a function defined on $\Uc$. If $H$ is smooth enough, we can associate with this function the Hamiltonian vector field $X_H(z)$ defined by
$$
X_H(z) = J \nabla H(z) 
$$
where $J$ is the symplectic operator on $\Pc_s$ induced  by the symplectic form \eqref{Esymp}.

For two functions $F$ and $G$, the Poisson Bracket is defined as
$$
\{F,G\} = \nabla F^T J \nabla G = i \sum_{a \in \Nc} \frac{\partial F}{\partial \eta_j}\frac{\partial G}{\partial \xi_j} -  \frac{\partial F}{\partial \xi_j}\frac{\partial G}{\partial \eta_j}.  
$$

We say that $z\in \Pc_s$ is {\em real} when $z_{\overline{j}} = \overline{z_j}$ for any $j\in \Zc$. In this case, $z=(\xi,\bar\xi)$ for some $\xi\in \C^{\Nc}$. Further we say that a Hamiltonian function $H$ is 
 {\em real } if $H(z)$ is real for all real $z$. 
 
\begin{definition}
Let $s \geq 0$, and let $\Uc$ be a neighborhood of the origin in $\Pc_s$. We denote by $\Hc^s(\Uc)$ the space of real Hamiltonian $H$ satisfying 
$$
H \in \mathcal{C}^{\infty}(\Uc,\C), \quad \mbox{and}\quad 
X_H \in \mathcal{C}^{\infty}(\Uc,\Pc_s). 
$$
\end{definition}

With a given function $H \in \Hc^s(\Uc)$, we associate the Hamiltonian system
$$
\dot z = J \nabla H(z)
$$
which can be written
\begin{equation}
\label{Eham2}
\left\{
\begin{array}{rcll}
\dot\xi_a &=& \displaystyle - i \frac{\partial H}{\partial \eta_a}(\xi,\eta) & a \in \Nc\\[2ex]
\dot\eta_a &=& \displaystyle i \frac{\partial H}{\partial \xi_a}(\xi,\eta)& a \in \Nc.  
\end{array}
\right.
\end{equation}
In this situation, we define the flow $\varphi_H^t(z)$ associated with the previous system (for times $t \geq 0$ depending on $z \in \Uc$). Note that if $z = (\xi,\bar \xi)$ and using the fact that $H$ is real valued, the flow $(\xi^t,\eta^t) = \varphi_H^t(z)$ satisfies for all time where it is defined the relation $\xi^t = \bar {\eta}^t$, where $\xi^t$ is solution of the equation 
\begin{equation}
\label{Eham1}
\dot\xi_a = - i \frac{\partial H}{\partial \eta_a}(\xi,\bar\xi), \quad a \in \Nc. 
\end{equation}
In this situation, introducing the real variables $p_a$ and $q_a$ such that
$$
\xi_a = \frac{1}{\sqrt{2}} (p_a + i q_a)\quad \mbox{and}\quad \bar{\xi}_a =  \frac{1}{\sqrt{2}} (p_a - i q_a),
$$
the system \eqref{Eham1} is equivalent to the system
$$
\left\{
\begin{array}{rcll}
\dot p_a &=& \displaystyle -  \frac{\partial H}{\partial q_a}(q,p) & a \in \Nc\\[2ex]
\dot q_a &=& \displaystyle  \frac{\partial H}{\partial p_a}(q,p),& 	a \in \Nc.  
\end{array}
\right.
$$
where $H(q,p) = H(\xi,\bar\xi)$. 

Note that the flow $\tau^t = \varphi_\chi^t$ of a real hamiltonian $\chi$ defines a symplectic map, i.e.  satisfies for all time $t$ and all point $z$ where it is defined
\begin{equation}
\label{Esympl}
(D_z  \tau^t)_z^T J (D_z\tau^t)_z = J
\end{equation}
where $D_z$ denotes the derivative with respect to the initial conditions. 

The following result is classic: 
\begin{lemma}
\label{Lchange}
Let $\Uc$ and $\Wc$ be two domains of $\Pc_s$, and let $\tau = \varphi_\chi^1 \in \mathcal{C}^{\infty}(\Uc,\Wc)$ be the flow of the real hamiltonian $\chi$. 
Then for $K \in \Hc^s(\Wc)$, we have 
$$
\forall\, z \in \Uc\quad
X_{K \circ \tau}(z) = 
(D_z\tau(z))^{-1}X_K(\tau(z)).
$$
Moreover, if $K$ is a real hamiltonian, $K \circ \tau$ is a real hamiltonian. 
\end{lemma}

\subsection{Function spaces}

We describe now the hypothesis needed on the Hamiltonian $H$. 

Let $\ell \geq 3$ be a given integer. For $\jb = (j_1,\ldots,j_r) \in \Zc^r$, we define $\mu(\jb)$ as the third largest integer between $|j_1|,\ldots,|j_r|$. Then we set $S(\jb) = |j_{i_r}| - | j_{i_{r-1}}| + \mu(\jb)$ where $|j_{i_r}|$ and $|j_{i_{r-1}}|$  denote the largest and the second largest integer between $|j_1|,\ldots,|j_r|$. 

Let $\ell \geq 3$. We consider $ \jb = (j_1,\ldots,j_\ell) \in \Zc^{\ell}$, and we set for all $i = 1,\ldots p$
$j_i = (a_i,\delta_i)$ where $a_i \in \Nc$ and $\delta_i \in \{\pm 1\}$. We define the moment $\Mc (\jb)$ of  the multi-index $\jb$ by
\begin{equation}
\label{EMcb}
\Mc(\jb) = a_{1} \delta_{1} + \cdots + a_\ell \delta_\ell.
\end{equation}
We then define the set of indices with zero moment
\begin{equation}
\label{EIr}
 \Ic_\ell =  \{  \jb = (j_1,\ldots,j_\ell) \in \Zc^{\ell}, \quad \mbox{with}\quad \Mc(\jb) = 0\}.
\end{equation}

In the following, for $\jb = (j_1,\ldots,j_\ell) \in \Ic_\ell$, we use the notation 
$$
z_\jb = z_{j_1}\cdots z_{j_\ell}. 
$$
Moreover, for $\jb = (j_1,\ldots,j_\ell) \in \Ic_\ell$ with $j_i = (a_i,\delta_i) \in \Nc \times\{ \pm 1\}$ for $i = 1,\ldots,\ell$, we set
$$
\overline \jb = (\overline{j}_1,\ldots,\overline j_\ell)\quad\mbox{with}\quad \overline{j}_i = (a_i,-\delta_i), \quad i = 1,\ldots,\ell.
$$
We recall  the following definition from \cite{Greb07}. 
\begin{definition}
Let $k \geq 3$, $M > 0$ and $\nu \in [0,+\infty)$, and let
$$
Q(z) = \sum_{\ell = 3}^k \sum_{\jb \in \Ic_\ell} Q_{\jb} z_{\jb}.
$$
We say that $Q \in \Tc_k^{M,\nu}$ if there exist a constant $C$ depending on $M$ such that 
\begin{equation}
\label{Ereg}
\forall\, \ell = 3,\ldots,k,\quad \forall\, \jb \in \Ic_\ell,\quad |Q_{\jb}| \leq C \frac{\mu(\jb)^{M+\nu}}{S(\jb)^M}. 
\end{equation}
\end{definition}
Note that $Q$ is a real hamiltonian if and only if
\begin{equation}
\label{Ereal}
\forall\, \ell = 3, \ldots,k, \quad \forall\, \jb \in \Ic_\ell, \quad Q_\jb = \overline{Q}_{\overline\jb}. 
\end{equation}

We have that $\Tc^{M,\nu}_k \in \Hc^s$ for $s \geq \nu +1/2$ (see  \cite{Greb07}). The best constant in the inequality \eqref{Ereg} defines a norm $\SNorm{Q}{\Tc^{M,\nu}_k}$ for which $\Tc^{M,\nu}_k$ is a Banach space. 
We set
$$
T_k^{\infty,\nu} = \bigcap_{M \in \N} \Tc^{M,\nu}_k. 
$$


\begin{definition}
A function $P$ is in the class $\Tc$ if
\begin{itemize}
\item $P$ is a real hamiltonian and exhibits a zero of order at least 3 at the origin. 
\item There exists $s_0 \geq 0$ such that for any $s \geq s_0$, $P \in \Hc^s(\Uc)$ for some neighborhood $\Uc$  of the origin in $\Pc_s$. 
\item For all $k \geq 1$, there exists $\nu \geq 0$ such that the Taylor expansion of degree $k$ of $P$ around the origin belongs to $\Tc_k^{\infty,\nu}$. 
\end{itemize}
\end{definition}
With previous notations, we consider in the following Hamiltonian functions of the form
\begin{equation}
\label{Edecomp}
H(z) = H_0(z) + P(z) = \sum_{a \in \Nc} \omega_a I_a(z)+ P(z),
\end{equation}
where for all $a\in \Nc$, 
$$
I_a(z) = \xi_a \eta_a
$$
are the {\em actions} associated with $a\in \Nc$ and
where $\omega_a \in \R$ are frequencies satisfying
\begin{equation}
\label{Eboundomega}
\forall\, a \in \Nc, \quad |\omega_a| \leq C |a|^m
\end{equation}
for some constants $C > 0$ and $m > 0$. 
The Hamiltonian system \eqref{Eham2} can hence be written
\begin{equation}
\label{Eham3}
\left\{
\begin{array}{rcll}
\dot\xi_a &=& \displaystyle - i \omega_a \xi_a - i \frac{\partial P}{\partial \eta_a}(\xi,\eta) & a \in \Nc\\[2ex]
\dot\eta_a &=& \displaystyle i \omega_a \eta_a + i \frac{\partial P}{\partial \xi_a}(\xi,\eta)& a \in \Nc.  
\end{array}
\right.
\end{equation}

\subsection{Rounded splitting methods}

When considering the numerical simulation of such hamiltonian system, many methods can be interpreted as splitting methods associated with   the decomposition \eqref{Edecomp}. This means that for small step size $h$, we approximate the flow $\varphi_H^h$ by the composed flow
$$
\varphi_H^h \simeq \varphi_{H_0}^h \circ \varphi_P^h. 
$$
For a given time $t$, and a small step size $h$ with $t = nh$, the approximation of $\varphi_H^t$ is then written
\begin{equation}
\label{Enumflot}
\varphi_H^t \simeq \left(\varphi_{H_0}^h \circ \varphi_P^h\right)^n. 
\end{equation}
We give examples of such schemes in the next section. 

In order to control the possible numerical instabilities due to the interaction of high frequencies, we introduce the following projection operator: Let $\eta > 0$ and $s$ be given, we define
$$
\Pi_{\eta,s} : \Pc_s \to \Pc_s
$$
by the formula
\begin{equation}
\label{Pi-1}
\forall\, j \in \Zc,\quad
\Big(\Pi_{\eta,s} z \Big)_{j} = \left\{
\begin{array}{rl}
z_j & \mbox{if}\quad |j|^s|z_j| \leq {\eta}\\[2ex]
0 & \mbox{if} \quad |j|^s|z_j| >  {\eta}.
\end{array}
\right.
\end{equation}

The goal of this paper is the studying of the long-time behavior of {\em rounded} splitting schemes associated with the operator
$$
\Pi_{\eta,s}\circ \varphi_{H_0}^h \circ \varphi_P^h
$$ 
to which we associate the numerical solution 
\begin{equation}
\label{Ezn}
z^n = \left(\Pi_{\eta,s} \circ \varphi_{H_0}^h \circ \varphi_P^h \right)^n(z^0). 
\end{equation}
Obviously, for $\eta = 0$, $\Pi_{\eta,s}$ is the identity operator. 

In the following, we show a normal form result on the abstract splitting method
$$
\varphi_{H_0}^h \circ \varphi_P^h
$$
and then draw some dynamical consequences for the discrete solution \eqref{Ezn}.

\section{Statement of the result and applications}

\subsection{Main result}

Let  $\jb = (j_1,\ldots,j_r)\in \Zc^r$, and denote by $j_i = (a_i,\delta_i) \in \Nc \times \{\pm 1\}$ for $i = 1,\ldots,r$. We set 
$$
\Omega(\jb) = 
\delta_1\omega_{a_1} + \cdots  + \delta_r\omega_{a_r}. 
$$
We say that $\jb = (j_1,\ldots,j_r) \in \Ic_r$ depends only of the action and we write $\jb \in \Ac_r$ if $r$ is even and if we can write 
$$
\forall\, i = 1,\ldots r/2,\quad
j_{i} = (a_i,1), \quad\mbox{and}\quad j_{i + r/2} = (a_i,-1)
$$
for some $a_i \in \Nc$. 
Note that in this situation, 
$$
\begin{array}{rcl}
z_\jb = z_{j_1}\cdots z_{j_r} &=& \xi_{a_1}\eta_{a_1} \cdots \xi_{a_{r/2}} \eta_{a_{r/2}}\\[2ex]
&=& I_{a_1} \cdots I_{a_{r/2}}
\end{array}
$$
where for all $a \in \Nc$, 
$$
I_{a}(z) = \xi_a \eta_a
$$
denote the action associated with the index $a$. Note that if $z$ satisfies the condition $z_{\overline{j}} = \overline{z_j}$ for all $j \in \Zc$, then we have $I_a(z) = |\xi_a|^2$. For odd $r$, $\Ac_r$ is the empty set.

We will assume now that the step size $h$ satisfies the following property:
\begin{hypothesis}\label{H1}
For all $r \in \N$, there exist constants $\gamma^*$ and $\alpha^*$ such that $\forall\, N \in \N^*$ and 
$\forall\jb = (j_1,\ldots,j_r) \notin  \Ac_r$, 
\begin{equation}
\label{nonres2}
|j_1|,\ldots,|j_r| \leq N\quad \Longrightarrow \quad |1 - e^{ih\Omega(\jb)}| \geq \frac{h \gamma^*}{N^{\alpha^*}}. 
\end{equation}
\end{hypothesis}

\begin{theorem}
\label{Tmain}
Assume that $P \in \Tc$ and $h < h_0$ satisfies the condition \eqref{nonres2}. Let $r \in \N^*$  be fixed. Then there exists a constant $s_0$ depending on $r$ such that for all $s > s_0$, there exist constants $C$ and $\varepsilon_0$ depending on $r$ and $s$ 
such that 
the following holds: For all $\varepsilon < \varepsilon_0$ and for all $z^0\in \Pc_{2s}$ real such that $\Pi_{\eta,s}z^0 = z^0$ with $\eta=\varepsilon^{r+1/4}$ and 
$$
\Norm{z^0}{s} \leq \varepsilon \quad \mbox{and} \quad \Norm{z^0}{2s} \leq 1,
$$ 
if we define 
\begin{equation}
\label{Eseuil}
z^n = \big(\Pi_{\eta,s}\circ\varphi_{H_0}^h \circ \varphi_P^h \big)^n (z^0)\quad \mbox{with}\quad \eta = \varepsilon^{r+1/4}
\end{equation}
then
we have $z^n$ is still real, 
and moreover
\begin{equation}
\label{Eresnorm}
\Norm{z^n}{s} \leq 2 \varepsilon \quad \mbox{for}\quad n \leq \frac{1}{\varepsilon^{r-2}},
\end{equation}
and  
\begin{equation}
\label{Eresact}
\sum_{a \in \Nc} |a|^{2s}| I_a(z^n) - I_a(z^0)| \leq \varepsilon^{5/2}\quad \mbox{for}\quad n \leq \frac{1}{\varepsilon^{r-2}}
\end{equation}
\end{theorem}
The proof is postponed to section \ref{proofTmain}.
\begin{remark}As $r$ is arbitrary, the  condition  $\eta = \varepsilon^{r+1/4}$ implies that $\Pi_{\eta,s}$ is $\varepsilon^{r+1/4}$ close to the identity in  $\Pc_{s}$ (cf. \eqref{Pi-1}). From the practical point of view, we can always assume that $\varepsilon^r$ is beyond the round-off error, so that we can consider that \eqref{Eseuil} coincides with the numerical solution associated with the splitting method. The full understanding of the {\em real} numerical phenomenon taking into account the round-off error is clearly out of the scope of this paper. 
\end{remark}
\begin{remark} The condition $\Norm{z_0}{2s}\leq1$ together with $\Pi_{\eta,s}z^0 = z^0$ implies that $z_j^0=0$ for $j$ large enough which is actually the assumption we need.
\end{remark}

\subsection{Verification of the non resonance condition}

We assume that the frequencies $\omega_a$, $a\in \Nc$ fulfill the following condition (in the next section we will verify this condition in different concrete case): 
\begin{hypothesis}\label{H2}
For all $r \in \N$, there exist  constants $\gamma(r)$ and $\alpha(r)$ such that $\forall\, N \in \N^*$ and 
$\forall\jb = (j_1,\ldots,j_r) \notin  \Ac_r$, 
\begin{equation}
\label{nonres}
|j_1|,\ldots,|j_r| \leq N\quad   \Longrightarrow |\Omega(\jb)| \geq \frac{ \gamma}{N^{\alpha}}. 
\end{equation}
\end{hypothesis}

The next result shows that under the previous hypothesis, the condition \eqref{nonres2} is generic. See \cite{Shan00,HLW} for similar statements. 

\begin{lemma}
Assume that Hypothesis \ref{H2} holds , and let $h_0$ and $r$ be  given numbers. Let $\gamma$ and $\alpha$ be such that \eqref{nonres} holds and assume that $\gamma^* \leq  (2/\pi) \gamma$, $\alpha^* \geq \alpha + m\sigma + r$ with $\sigma> 1$ and $m$ the constant appearing in \eqref{Eboundomega}, then we have 
$$
\mbox{\rm meas}\{\, h < h_0\, | \, h \mbox{ does not satisfy }  \eqref{nonres2}\, \} \leq C \frac{\gamma^*}{\gamma} h_0^{1+\sigma}
$$
where $C$ depends on $\sigma$ and  $r$. As a consequence the set  $$
Z(h_0) = \{\, h < h_0\, | \, h \mbox{ satisfies Hypothesis }  \ref{H1}\, \}
$$
is a dense open subset of $(0,h_0)$.

\end{lemma}
\begin{Proof}
Denote $$\Rc(h_0, \gamma^*,\alpha^*)= \{\, h < h_0\, | \, h \mbox{ does not satisfy }  \eqref{nonres2}\}.$$
Assume that $h \in \Rc(h_0, \gamma^*,\alpha^*)$. There exist $N > 1$ and $\jb\notin\Ac_r$ such that 
$$
|j_1|,\ldots,|j_r| \leq N\quad \mbox{and}\quad | 1 - e^{ih\Omega(\jb)}| <  \frac{h \gamma^*}{N^{\alpha^*}}.
$$
For this $\jb$, there exist an $\ell \in \Z$ such that 
\begin{equation}
\label{E13}
| 1 - e^{ih\Omega(\jb)}| \geq \frac{2}{\pi} | 2\pi \ell - h\Omega(\jb)|. 
\end{equation}
If $\ell = 0$, the previous inequality and \eqref{nonres} imply 
$$
| 1 - e^{ih\Omega(\jb)}| \leq \frac{2}{\pi} h \frac{\gamma}{N^\alpha}
$$
which is impossible with the assumptions on $\gamma^*$ and $\alpha^*$. Hence, we can assume $\ell\neq 0$. 
Eqn. \eqref{E13} implies 
$$
\frac{2|\Omega(\jb)|}{\pi} \Big| \frac{2\pi \ell}{\Omega(\jb)} - h \Big| < \frac{h \gamma^*}{N^{\alpha^*}}
$$
and using \eqref{nonres}
$$
\Big| \frac{2\pi \ell}{\Omega(\jb)} - h | \leq \frac{h\pi\gamma^*}{2\gamma} \frac{1}{N^{\alpha^* - \alpha}}.  
$$
Moreover, we have for this $\ell$
$$
| 2\pi \ell - h\Omega(\jb)| \leq \pi
$$
whence using \eqref{Eboundomega}
$$
2\pi | \ell| \leq \pi + C h_0 N^m
$$
where $C$ is a constant depending on $r$. This implies 
$$
| \ell| - \frac12 \leq \frac{C}{2\pi} h_0 N^m. 
$$
Hence, $\Rc(h_0, \gamma^*,\alpha^*)$ is included in the union of balls of center 
$$
\frac{2\pi\ell}{\Omega(\jb)}, \quad\mbox{with} \quad |j_1|,\ldots,|j_r| \leq N,\quad |\ell| \leq \pi + Ch_0 N^m, \quad \ell \neq 0
$$
and radius
$$
\frac{h_0\pi\gamma^*}{2\gamma} \frac{1}{N^{\alpha^* - \alpha}}
$$
Hence, we have for $\sigma > 1$
\begin{equation*}
\begin{split}
\mbox{\rm meas}(\Rc(h_0, \gamma^*,\alpha^*)) &\leq \sum_{|j_i|\leq N} \, \sum_{|\ell| - \frac12 \leq \frac{C}{2\pi} h_0 N^m} \frac{h_0\pi \gamma^*}{2 \gamma} \frac{1}{N^{\alpha^* - \alpha}}\\[1ex]
&\leq \sum_{|j_i|\leq N} \,\sum_{\ell \in \Z^*}
\Big(\frac{1}{|\ell| - \frac12}\Big)^{\sigma} \frac{h_0\pi \gamma^*}{2 \gamma}  \frac{1}{N^{\alpha^* - \alpha - m\sigma}} \Big(
\frac{Ch_0}{2\pi}
\Big)^\sigma. \\[1ex]
&\leq C \frac{\gamma^*}{\gamma} h_0^{1 + \sigma} \frac{1}{N^{\alpha^* - \alpha - m\sigma - r}}.
\end{split}
\end{equation*}
Furthermore
$$ \mbox{\rm meas}(\cap_{\gamma^*>0}\Rc(h_0, \gamma^*,\alpha^*))=0$$
and thus $Z(h_0)$ has full measure.
\end{Proof}

\subsection{Examples}
In this section we  present two examples, other examples like the Klein Gordon equation on the sphere (in the spirit of \cite{BDGS}) or the nonlinear Schr\"odinger operator with harmonic potential (in the spirit of \cite{GIP}) could also be considered with these technics but will require some additional efforts.

\subsubsection{Schr\"odinger equation on the torus}

We first consider  non linear Schr\"odinger equations of the form
\begin{equation}
\label{E1}
i \partial_t \psi = - \Delta \psi + V \star\psi + \partial_2g(\psi,\bar \psi),\quad x \in \T^d
\end{equation}
where $V\in C^\infty(\T^d,\R)$,  $g\in C^\infty(\Uc,\C)$ where $\Uc$ is a neighborhood of the origin in $\C^2$. We assume that $g(z,\bar z) \in \R$, and that $g(z,\bar z) = \mathcal{O}(|z|^3)$. 
The corresponding hamiltonian functional is given by 
$$
H(\psi,\bar\psi) = \int_{\T^d} | \nabla \psi | ^2 + \bar\psi (V \star \psi) + g(\psi,\bar\psi) \, \dd x
$$

Let $\phi_{a}(x) = e^{i a\cdot x}$, $a \in \Z^d$ be the Fourier basis on $L^2(\T^d)$. With the notation 
$$
\psi = \Big(\frac{1}{2\pi}\Big)^{d/2}\sum_{a\in \Z^d} \xi_{a} \phi_{a}(x) \quad \mbox{and}\quad
\bar \psi = \Big(\frac{1}{2\pi}\Big)^{d/2}\sum_{a\in \Z^d} \eta_{a}\bar \phi_{a}(x)
$$
 the hamiltonian associated with the equation \eqref{E1} can be (formally) written
\begin{equation}
\label{E2}
H(\xi,\eta) = 
\sum_{a \in \Z^d} \omega_a \xi_a \eta_a  + \sum_{r \geq 3}\, 
\sum_{\ab,\bb } P_{\ab\bb}\,  \xi_{a_1} \cdots \xi_{a_p}\eta_{b_1}\cdots \eta_{b_q}.
\end{equation}
Here  $\omega_{a} =|a|^2 +\hat V_a$ satisfying \eqref{Eboundomega} with $m = 2$ and are the eigenvalues of the operator 
$$
\psi\mapsto  - \Delta\psi + V \star\psi.
$$
Note that in \eqref{E2} the sum is made over the set of multi-indices 
\begin{multline*}
\{  (\ab,\bb) = (a_1,\ldots,a_p,b_1,\ldots,b_q) \in (\Z^d)^{p}\times (\Z^d)^q\, \quad \mbox{with}\quad p + q = r\\[1ex] 
 \mbox{and}\quad a_1 + \cdots + a_p - b_1 - \ldots - b_q = 0\},
\end{multline*}
which corresponds to the set \eqref{EIr} in variables $(z_j)_{j \in \Zc}$ (here we set $\Nc = \Z^d$). 

The relation $H(\xi,\bar\xi) \in \R$ is equivalent to the fact that the coefficients $P_{\ab\bb}$ satisfy $P_{\ab\bb} = \overline{P}_{\bb\ab}$ which correspond to the hypothesis \eqref{Ereal} in variables $(z_j)_{j \in \Zc}$. The fact that the nonlinearity 
$P$ belongs to $\Tc$ can be verified using the regularity of $g$ and the properties of the basis functions $\phi_a$, see \cite{Greb07,BG06}. In this situation, it can be shown  that the Hypothesis \ref{H1} is fulfilled for a large set of potential $V$ (see \cite{BG06} or \cite{Greb07}). 

The numerical implementation of splitting method is very easy in the case of Eqn. \eqref{E1}:  The part corresponding to the equation
$$
i\partial_t \psi = -\Delta \psi + V \star \psi
$$
is easily solved in terms of Fourier coefficients, while the non linear part 
$$
i\partial_t \psi = \partial_2g(\psi,\bar\psi)
$$
is a simple differential equation with fixed $x \in \T^d$. The use of a fast Fourier transform allows to compute alternatively the solution of the linear part and the solution of the non-linear part.

\subsubsection{Wave equation on the circle}

We consider the wave equation on the circle
$$
u_{tt} - u_{xx} + m u = g(u), \quad x \in \T^1, \quad t \in \R, 
$$
where $m$ is a non negative real constant and $g$ a smooth real valued  function.
Introducing the variable $ v = u_t$, the corresponding hamiltonian can be written
$$
H(u,v) =\int_{\T} \frac12 (v^2 +  u_x^2 +  m u^2) + G(u) \, \dd x,
$$
where $G$ is such that $\partial_u G = g$. Let $A := (-\partial_{xx} + m)^{1/2}$, and define the variables $(p,q)$ by
$$
q := A^{1/2} u, \quad\mbox{and}\quad p = A^{-1/2}v. 
$$
Then the Hamiltonian can be written 
$$
H = \frac12 \big(\langle Ap,q\rangle_{L^2} + \langle Aq,q\rangle_{L^2} \big) + \int_{\T} G(A^{-{1/2}}q) \, \dd x. 
$$
Let $\omega_a=\sqrt{|a|^2+m}$, $a \in \N =: \Nc$ be the eigenvalues of the operator $A$, and $\phi_a$ the associated eigenfunctions. Plugging the decompositions 
$$
q(x) = \sum_{a \in \N} q_a \phi_a(x) \quad \mbox{and} \quad p(x) = \sum_{a\in \N} p_a \phi_a(x)
$$
into the hamiltonian functional, we see that it takes the form 
$$
H = \sum_{a \in \N}\omega_a \frac{p_a^2 + q_a^2}{2} + P
$$
where $P$ is a function of the variables $p_a$ and $q_a$. 
Using the complex coordinates 
$$
\xi_a = \frac{1}{\sqrt{2}} ( q_a + ip_a) \quad\mbox{and}\quad 
\eta_a = \frac{1}{\sqrt{2}} ( q_a - ip_a)
$$
the hamiltonian function can be written under the form \eqref{Edecomp} with $P \in \Tc$ (see \cite{Bam07, Greb07}). As in the previous case, it can be shown  that the condition \eqref{nonres} is fulfilled for a l set of constant $m$ of full measure (see \cite{BG06, Bam07}).  

In this situation, the symmetric Strang splitting scheme 
$$
\varphi_{P}^{h/2} \circ \varphi_{H_0}^h \circ \varphi_{P}^{h/2}
$$
corresponds to the Deuflhard's method \cite{Deuf79}. Considering now the Hamiltonian 
$$
H(z) = H_0(z) + P(\Phi(h H_0) z)
$$
where $\Phi(x)$  a smooth function that is real, bounded, even and such that $\Phi(0) = 1$, and where 
$$
\Phi(h H_0)(\xi,\eta)= \Big(\Phi(h \omega_a) (\xi_a,\eta_a)\Big)_{a \in \Nc}, 
$$
then 
the  splitting schemes associated with this decomposition coincide with the (abstract) symplectic mollified impulse methods (see \cite[Chap. XIII]{HLW} and \cite{CHL08c}). The fact that  $\Phi$ is bounded makes that the functional $z \mapsto P(\Phi(H_0) z)$ obviously belongs to $\Tc$.

\section{A normal form result}

\subsection{Normal form}

\begin{definition}
Let $N > 0$ be a real number. For a given multi-index $\jb \in \Zc^r$, let $i_p$ be the permutation such that 
$$
|j_{i_1}| \leq \cdots \leq |j_{i_r}|. 
$$
We define the set 
$$
\Jc_r(N) = \{ \jb \in \Ic_r\, \, |j_{i_r}| \leq (r-1) N \quad \mbox{and} \quad |j_{i_{r-1}}| \leq N\}. 
$$
\end{definition}
\begin{lemma}
\label{L1}
Let $r \geq 3$, and 
assume that $\jb \notin \Jc_r(N)$. Then $\jb$ contains at least two indices with modulus greater than $N$. 
\end{lemma} 
\begin{Proof} 
Let $\jb = (j_1,\ldots,j_r) \in \Ic_r\backslash \Jc_r(N)$. We have
$\Mc(\jb) = 0$
where $\Mc(\jb)$ is defined in \eqref{EMcb}. 
Assume that there exists only one index of modulus greater than $N$. We can assume that $|j_1| > N$ and hence all the other indices are of modulus $\leq N$ (in particular, with the previous notation, we have $j_1 = j_{i_r}$). Hence we have
$$
|j_1| \leq |j_2| + \cdots + |j_r|  \leq  (r-1)N 
$$
and this implies that $\jb \in \Jc_r(N)$ which is a contradiction. 
\end{Proof}

We motivate now the definition of normal form terms we introduce in the sequel. For a given number $N$ and $z \in \Pc_s$ we define 
$$
\Bc_{s}^{N}(z) = \sum_{|a| \leq N} |a|^{2s} \xi_a \eta_a
$$
and
$$
\Rc_{s}^{N}(z) = \sum_{|a| > N} |a|^{2s} \xi_a \eta_a
$$
so that 
$$
\Norm{z}{s}^2 = \Bc_{s}^{N}(z) + \Rc_{s}^{N}(z). 
$$
\begin{proposition}
\label{Pcrux}
Let $N \in \N$ and $r \geq 3$. 
 Assume that the homogeneous polynomial
$$
Z= \sum_{\jb \in \Ic_r\backslash\Jc_r(N)} Z_{\jb} z_{\jb} 
$$
defines an element of $\Tc_r^{M,\nu}$ for some constants $M$ and $\nu$. 
Then we have for all $s > 2\nu + 4$, $M > s +2$ and for all $z \in \Pc_s(\C)$,
\begin{equation}
 \label{ENsP}
|\{\Bc_s^{N},Z\}(z)| \leq  C_0 \SNorm{Z}{\Tc_r^{M,\nu}} N^{\nu + 2 + d/2- s}\Norm{z}{s}^{r-2} \Rc_s^N(z). 
\end{equation}
and 
\begin{equation}
\label{ENsP2}
\forall\, a \in \Nc,\quad |a|\leq N, \quad|\{I_a,Z\}| \leq  C_0 \SNorm{Z}{\Tc_r^{M,\nu}} N^{\nu + 2 - s}\Norm{z}{s}^{r-2} \Rc_s^N(z). 
\end{equation}
Moreover
\begin{equation}
\label{ERsP}
|\{\Rc_s^{N},Z\}(z)| \leq  C_0 \SNorm{Z}{\Tc_r^{M,\nu}} \Norm{z}{s}^{r-2} \Rc_s^N(z)
\end{equation}
where $C_0$ is a constant depending on $s$, $r$ and the dimension $d$ of $\Nc = \N^d$ of $\Z^d$.
\end{proposition}

The proof of this proposition is given in the Appendix. 

\begin{definition}
An element $Z \in \Tc_{r}^{M,\nu}$ is said to be in normal form if we can write it
$$
Z = \sum_{\ell = 3}^{r} \sum_{\{\jb \in \Ac_\ell \cup \Ic_\ell \backslash \Jc_\ell(N) \}} Z_\jb z_{\jb}. 
$$
\end{definition}

In other words, a normal form term either depends only on the actions or contains at least two terms with index greater than $N$ (cf. Lemma \ref{L1}). 


\subsection{Statement of the normal form result}
In the following, we set
$$
B_s(\rho) = \{ z \in \Pc_s \, | \, \Norm{z}{s} \leq \rho\}.
$$

\begin{theorem}
\label{TT}
Assume that $P \in \Tc$ and $h < h_0$ satisfies the Hypothesis  \ref{H1}. Let $r_0 \geq 3$ be fixed. Then there exist constants $s_0$, $\beta$ and $N_0$ such that for all $s \geq s_0$, there exists constants  $C$  and for all $N \geq N_0$ there exists a canonical transformation $\tau$ from $B_s({\rho})$ into $B_s({2\rho})$ with $\rho = (C N)^{-\beta}$ satisfying for all $z \in B_s({\rho})$, 
\begin{equation}
\label{Eestbasis}
\Norm{\tau(z) - z}{s}\leq (CN)^{\beta}\Norm{z}{s}^2 \quad\mbox{and}\quad \Norm{\tau^{-1}(z) - z}{s}\leq (CN)^{\beta}\Norm{z}{s}^2
\end{equation}
and such that the restriction of $\tau$ to the high modes is the identity, i.e.
\begin{equation}
\label{Ehighmodes}
(\tau(z))_j = z_j \quad \mbox{for}\quad |j| > (r_0-1)N.
\end{equation}
Moreover, $\tau$ puts $\varphi_H^h$ in normal form up to order $r_0$ in the sense that
\begin{equation}
\label{Efactori}
\varphi_{H_0}^h \circ \varphi_P^h \circ \tau = \tau \circ \varphi_{H_0}^h \circ \psi
\end{equation}
where $\psi$ is the solution at time $\lambda = 1$ of a non-autonomous hamiltonian $h Z(\lambda) + R(\lambda)$ with
\begin{itemize}
\item $Z(\lambda) \in \mathcal{C}([0,1],\Tc_{r_0}^{M_1,\nu_1})$ for some $M_1$ and $\nu_1$ depending on $P$, $r_0$, $s$ and $h_0$, and for all $\lambda \in [0,1]$, $Z(\lambda)$ is a real polynomial of degree $r$ under normal form such that
\begin{equation}\label{ZZ}
 |Z(\lambda)|_{T_{r_0}^{M_1,\nu_1}} \leq (C N)^{\beta}.
\end{equation}
 
\item $R(\lambda) \in \mathcal{C}([0,1], \Hc^s(B_s({\rho})))$ with $\rho \leq (C N)^{-\beta}$ has a zero of order $r_0+1$ at the origin and satisfies and for all $z \in B_s(\rho)$, 
\begin{equation}\label{RR}
\forall\, \lambda \in [0,1], \quad\Norm{X_{R(\lambda)}(z)}{s} \leq (C N)^{\beta}\Norm{z}{s}^{r}. 
\end{equation}
\end{itemize}
\end{theorem}
The proof is postponed to section \ref{5.4} and \ref{5.6}. We first verify that this normal form theorem has the dynamical consequences announced in Theorem \ref{Tmain}.

\subsection{Proof of the main Theorem \ref{Tmain}}\label{proofTmain}

We now give the proof of Theorem \ref{Tmain}.

First, let us note that as the Hamiltonian functions $H_0$, $P$, $Z$ and $R$ are real hamiltonians, and by definition of $\Pi_{\eta,s}$ (which is symmetric in $\xi$ and $\eta$), it is clear that there exist $\xi^n \in \C^\Nc$ such that for all $n$, we have $z^n = (\xi^n,\bar\xi^n)$. 

Let $r_0 = r$, and 
let $C = C(r,s)$, $s_0 = s_0(r)$ and $\beta = \beta(r)$ be the constants appearing in Theorem \ref{TT}. We can always assume that  
\begin{equation}
\label{Econd}
s_0(r) \geq 2 (r+1) \beta(r). 
\end{equation}
Let $s \geq s_0(r)$. 
Let $\varepsilon_0$ be such that
\begin{equation}
\label{Eepsilon0}
\varepsilon_0^{1/2} \leq \frac{(r-1)^{s} }{2 C^s}. 
\end{equation}
Let $\varepsilon < \varepsilon_0$. We define $N$ such that
\begin{equation}
\label{EchoixN}
(CN)^{\beta} = \varepsilon^{-1/2}. 
\end{equation}
Notice that the assumption $\Norm{z^0}{s}\leq \varepsilon$ implies $z^0\in B_s(\rho)$ with $\rho=(CN)^{-\beta}$.
Furthermore  we have $\Norm{z^0}{2s}\leq 1$ and together with $\Pi_{\eta,s} z^0 = z^0$, this hypothesis implies that $z^0_j=0$ for $j$ large enough. 
Actually let $j \in \Zc$ be such that $|j| > (r-1)N$, we have
$$
|j^s||z^0_j|\leq |j|^{-s}\leq ((r-1)N)^{-s}.
$$
using \eqref{EchoixN} we have $N = C^{-1} \varepsilon^{-\frac{1}{2\beta}}$ and hence
$$
|j^s||z^0_j| \leq (r-1)^{-s} C^s \varepsilon^{\frac{s}{2\beta}}.  
$$
Now condition \eqref{Econd} implies that  $\frac{s}{2\beta} \geq r +1$ and hence (as we can always assume that $\varepsilon < 1$), 
$$
\varepsilon^{\frac{s}{2\beta}}\leq \varepsilon^{r +1}.
$$
Therefore we get using \eqref{Eepsilon0}
$$
|j^s||z^0_j| \leq \varepsilon^{r+1/2} \Big( (r-1)^{-s} C^s \varepsilon^{1/2} \Big) \leq \varepsilon^{r+1/2}.
$$
As $\Pi_{\eta,s} z^0 = z^0$ and $\eta=\varepsilon^{r+1/4}$, this implies that 
$$
\forall\, |j| > (r-1)N,\quad z_j^0 = 0. 
$$

Let $\tau$ defined by Theorem \ref{TT}, and let $y^n = \tau^{-1}(z^n)$. As $\tau$ is the flow of a real hamltonian, there exist $\zeta^n \in \C^\Nc$ such that $y^n = (\zeta^n,\bar\zeta^n)$ for all $n$.  By definition, we have
\begin{equation}\label{star}
\forall\, n \geq 0,\quad 
y^{n+1} = \big(\tau^{-1} \circ \Pi_{\eta,s} \circ \tau\big) \circ \big(\varphi_{H_0}^h \circ \psi\big) (y^n). 
\end{equation}
and as $\tau$ is the identify for high modes (see \eqref{Ehighmodes}), we have $y_j^0 = 0$ for $|j| > (r-1)N$.

Using the definition of $N$, the transformation $\tau$ in the previous Theorem satisfies (taking $\rho := 2\varepsilon < \sqrt{\varepsilon}$): for all $z$ such that $\Norm{z}{s} \leq  2\varepsilon$, 
\begin{equation}
\label{Etransfo}
\begin{array}{rcl}
\Norm{\tau^{-1}(z) - z}{s} & \leq &  \varepsilon^{-1/2} \Norm{z}{s}^2 \\[1ex]
&\leq& 4 \varepsilon^{3/2}\\[1ex]
&\leq& \frac14 \varepsilon
\end{array}
\end{equation}
provided $\varepsilon_0$ is sufficiently small. 
Hence, we have $\Norm{y^0}{s} = \Norm{\tau^{-1}(z^0)}{s} \leq \frac54\varepsilon$. 

We will show by induction that the following holds for all $n \in \N$: 
\begin{itemize}
\item[(i)] $\Norm{y^n}{s}^2 \leq \Norm{y^0}{s}^2 + 2n \varepsilon^{r+1/8}$
\item[(ii)] $y^n_j = 0$ for $j \geq (r-1)N$. 
\end{itemize}

These assumptions are satisfied for $n = 0$. Assume that they hold for $n \geq 0$. 

Let $\psi$ be the application defined by Theorem \ref{TT} and $\psi^\lambda(z)$ be the flow associated with the Hamiltonian $h Z(\lambda) + R(\lambda)$ defining the application $\psi$ for $\lambda = 1$.

Using the results of Lemmas \ref{LP1} and \ref{LP2} below, we easily see that there exists a constant $c$ depending on $r$ such that for all $\lambda \in [0,1]$, $\Norm{\psi^{\lambda}(y^n ) }{s} \leq c \varepsilon$. 

Let $N_1 = (r - 1)N$. We have by hypothesis that 
$\Rc_s^{N_1}(y^n) = 0$. Furthermore
$$
\frac{\dd}{\dd \lambda}\Rc_s^{N_1}(\psi^\lambda(y^n)) =\{\Rc_s^{N_1},hZ+R\}.
$$
Thus using the equation \eqref{ERsP}, \eqref{ZZ} and \eqref{RR}  we get 
\begin{multline*}
\left|\frac{\dd}{\dd \lambda}\Rc_s^{N_1}(\psi^\lambda(y^n)) \right| \\[1ex]
\leq 
C_1 N^{\beta} \Rc_s^{N_1}(\psi^\lambda(y^n)) \big(\Norm{\psi^\lambda(y^n)}{s} + \Norm{\psi^\lambda(y^n)}{s}^{r-2}\big) + C_1 N^{\beta} \Norm{\psi^{\lambda}(z)}{s}^{r+1} 
\end{multline*}
for some constant $C_1$ depending on $r$ and $s$. Hence we have
$$
\left|\frac{\dd}{\dd \lambda}\Rc_s^{N_1}(\psi^\lambda(y^n)) \right| \leq 
C_1 \big( \varepsilon^{1/2}\Rc_s^{N_1}(\psi^\lambda(y^n))  + \varepsilon^{r + 1/2}\big).
$$
where $C_1$ depends on $r$ and $s$. Using the Gronwall Lemma, we obtain for all $\lambda \in [0,1]$
$$
\Rc_s^{N_1}(\psi^\lambda(y^n))\leq \varepsilon^{r + 1/2}  C_1 e^{\varepsilon^{1/2}\lambda C_1}. 
$$
We can always assume that $\varepsilon_0^{1/4} C_1 e^{\varepsilon_0^{1/2} C_1} < 1$. Hence we get 
$$
\forall\, \lambda \in [0,1], \quad \Rc_s^{N_1}(\psi^\lambda(y^n))\leq \varepsilon^{r+1/4}.  
$$

On the other hand, using \eqref{ENsP} we have  
\begin{multline*}
\left|\frac{\dd}{\dd \lambda}\Bc_s^{N_1}(\psi^\lambda(y^n)) \right| \leq\\[1ex] 
C_1 N^{\beta + \nu + 2 + d/2 - s} \Rc_s^{N_1}(\psi^\lambda(y^n)) \big(\Norm{\psi^\lambda(z)}{s} + \Norm{\psi^\lambda(z)}{s}^{r-2}\big) + C_1N^{\beta} \Norm{\psi^{\lambda}(z)}{s}^{r+1}
\end{multline*}
for some constant $C_1$ depending on $r$ and $s$. 
We can always assume that $s > \beta + \nu + d/2 + 2$. Using the previous estimates, we get
$$
\left|\frac{\dd}{\dd \lambda}\Bc_s^{N_1}(\psi^\lambda(y^n)) \right| 
\leq C_1 \varepsilon^{r + 1/2}
$$ and then
$$\Bc_s^{N_1}(\psi^\lambda(y^n))\leq  \Norm{y^n}{s}^2 + C_1 \varepsilon^{r + 1/2}.$$
Let $\tilde{y}^{n} = \varphi_{H_0}^h \circ \psi(y^n)$. 
As for all $z$ we have $\Norm{z}{s}^2 = \Bc_s^{N_1}(z) + \Rc_s^{N_1}(z)$ and as $\varphi_{H_0}^h$  preserves all the actions, therefore
\begin{equation}\label{xx}
\Norm{\tilde{y}^{n}}{s}^2 \leq \Norm{y^n}{s}^2 + C_1 \varepsilon^{r+1/2}\quad \mbox{and} \quad
\Rc_s^{N_1}(\tilde{y}^{n}) \leq \varepsilon^{r +1/4}.
\end{equation}
Now by construction (cf. \eqref{star})
$$
y^{n+1} = \tau^{-1} \circ \Pi_{\eta,s} \circ \tau(\tilde y^n). 
$$
As $\tau$ is the identity for modes $|j| > N_1$, we have 
$$
\Rc^{N_1}_s(\tau(\tilde{y}^n)) = \Rc^{N_1}_s(\tilde{y}^n) \leq \varepsilon^{r+1/4} = \eta.  
$$
Hence by definition of the projection $\Pi_{\eta,s}$ we  get that 
$$
\Big(\Pi_{\eta,s} \circ \tau(\tilde{y}^{n})\Big)_j = 0,\quad \mbox{for}\quad |j| > (r-1)N=N_1.
$$
As $\tau^{-1}$ is the identity for modes greater than $(r-1)N$, this shows (ii) for $n+1$, i.e.
$$
y^{n+1}_j = 0 \quad \mbox{for}\quad |j| > (r-1)N.
$$
Let $z$ be such that $z_{j} = 0$ for $|j| > (r-1)N$. We have 
$$
\Norm{\Pi_{\eta,s}z - z}{s} \leq \sum_{|j| \leq N} \eta \leq \eta N^d\leq \varepsilon^{r+1/4-d/2\beta}\leq \varepsilon^{r+1/8}
$$ 
since we can always assume $\beta>4d$.\\
Writing
$$
\tau^{-1} \circ \Pi_{\eta,s} \circ \tau = I + \tau^{-1} \circ (\Pi_{\eta,s} - I) \circ \tau
$$
and as $\tau$ leaves the set $(z_j)_{|j| \geq N_1}$ invariant, we get using \eqref{Eestbasis}, 
\begin{align*}
\Norm{y^{n+1}}{s}^2= \Bc_s^{N_1}(y^{n+1}) &\leq  \Bc_s^{N_1}(\tilde{y}^n) + ( 1 + (CN)^\beta\varepsilon^2) \eta N^d\\
& \leq \Bc_s^{N_1}(\tilde{y}^n) + \frac 3 2\varepsilon^{(r+1/8)}
\end{align*}
Thus we get
$$
\Norm{y^{n+1}}{s}^2 \leq \Norm{y^n}{s}^2 + c \varepsilon^{r+1/4} +\frac 3 2 \varepsilon^{r+1/8}
\leq \Norm{y^n}{s}^2 + 2\varepsilon^{r+1/8}.
$$

This shows (i) for $n +1$. 

In particular, for all $n \leq \varepsilon^{-r+2}$ we have (recall $\Norm{y^0}{s}\leq \frac54 \varepsilon$)
$$
\Norm{y^n}{s}^2 \leq \big(\textstyle\frac{5}{4}\varepsilon\big)^2 +  2\varepsilon^{2+1/8}
$$
and hence (provided $\varepsilon_0$ is small enough)
$$
\Norm{y^n}{s} \leq \textstyle\frac{7}{4}\varepsilon. 
$$
Now using \eqref{Etransfo} for the application $\tau$, we easily see that $\Norm{z^{n}}{} \leq 2 \varepsilon$ as long as  $n \leq \varepsilon^{r-2}$. This proves \eqref{Eresnorm}. 

The proof of \eqref{Eresact} is obtained similarly  using  \eqref{ENsP2} and we do not give the details here. 

\subsection{Formal equations}\label{5.4}

We give now the general strategy of the proof of the normal form Theorem \ref{TT},  showing in particular the need of working with non autonomous Hamiltonians and of considering the non resonance condition \eqref{nonres2}. 

We consider  a fixed step size $h$ satisfying \eqref{nonres2}. 
We consider the propagator
$$
\varphi_{H_0}^h \circ \varphi_{P}^h = \varphi_{H_0}^h \circ \varphi_{hP}^1. 
$$
We embed this application into the family of applications
$$
\varphi_{H_0}^h \circ \varphi_{hP}^\lambda, \quad \lambda \in [0,1]. 
$$
Formally, we would like to find a real hamiltonian $\chi = \chi(\lambda)$ and a real hamiltonian under  normal form $Z = Z(\lambda)$ and such that 
\begin{equation}
\label{eq:flots}
\forall\, \lambda \in [0,1]\quad \varphi_{H_0}^h \circ \varphi_{hP}^\lambda \circ \varphi_{\chi(\lambda)}^\lambda =  \varphi_{\chi(\lambda)}^\lambda \circ\varphi_{H_0}^h \circ \varphi_{hZ(\lambda)}^\lambda. 
\end{equation}
Let $z^0 \in \Pc_s(\C)$ and $z^\lambda = \varphi_{H_0}^h \circ \varphi_{hP}^\lambda \circ \varphi_{\chi(\lambda)}^\lambda (z^0)$. Deriving the previous equation with respect to $\lambda$ yields
\begin{multline*}
\frac{\dd z^\lambda}{ \dd \lambda} = 
(D_z \varphi_{H_0}^h)_{\varphi_{H_0}^{-h}(z^\lambda)} X_{hP}( \varphi_{H_0}^{-h}(z^\lambda) ) + \\[1ex]
(D_z( \varphi_{H_0}^h \circ \varphi_{hP}^\lambda ))_{ \varphi_{hP}^{-\lambda} \circ \varphi_{H_0}^{-h}(z^\lambda)}
X_{\chi } (\varphi_{hP}^{-\lambda} \circ \varphi_{H_0}^{-h}(z^\lambda)) . 
\end{multline*}
Using Lemma \ref{Lchange} that remains obviously valid for non autonomous hamiltonian, we thus have 
$$
\frac{\dd z^\lambda}{ \dd \lambda} = X_{A(\lambda)}(z^\lambda)
$$
where $A(\lambda)$ it the time dependent real hamiltonian given by 
$$
A(\lambda) = hP \circ \varphi_{H_0}^{-h}+ \chi(\lambda) \circ \varphi_{hP}^{-\lambda} \circ \varphi_{H_0}^{-h}. 
$$
Using the same calculations for the right-hand side, \eqref{eq:flots} is formally equivalent to the following equation (up to an integration constant)
\begin{equation}
\label{eq:tg}
\forall\, \lambda \in [0,1]\quad 
hP \circ \varphi_{H_0}^{-h} + \chi \circ \varphi_{hP}^{-\lambda} \circ \varphi_{H_0}^{-h} = \chi(\lambda) + hZ(\lambda) \circ \varphi^{-\lambda}_{\chi(\lambda)} \circ \varphi_{H_0}^{-h}. 
\end{equation}
which is equivalent to 
\begin{equation}
\label{eq:tg1}
\forall\, \lambda \in [0,1]\quad 
  \chi(\lambda) \circ \varphi_{H_0}^h - \chi(\lambda) \circ \varphi_{hP}^{-\lambda}= hP - hZ(\lambda) \circ \varphi^{-\lambda}_{\chi(\lambda)}   . 
\end{equation}
In the following, we will solve this equation in $\chi(\lambda)$ and $Z(\lambda)$ with a remainder term of order $r$ in $z$. 
So instead of \eqref{eq:tg1}, we will solve  in section \ref{5.6} the equation
\begin{equation}
\label{eq:tg2}
\forall\, \lambda \in [0,1]\quad 
  \chi(\lambda) \circ \varphi_{H_0}^h - \chi(\lambda) \circ \varphi_{hP}^{-\lambda}= hP - (h Z(\lambda) + R(\lambda)) \circ \varphi^{-\lambda}_{\chi(\lambda)} . 
\end{equation}
where the unknown are $\chi(\lambda)$, and $Z(\lambda)$ are polynomials of order $r$, with $Z$ under normal form, and where $R(\lambda)$ possesses a zero of order $r+1$ at the origin.

In the following, we formally write 
$$
\chi(\lambda) = \sum_{\ell = 3}^r \chi_{[\ell]}(\lambda)  :=  \sum_{\ell = 3}^r\sum_{\jb \in \Ic_\ell} \chi_\jb(\lambda)  z_{\jb}
$$
and 
$$
Z(\lambda) = \sum_{\ell = 3}^r Z_{[\ell]}(\lambda)  :=  \sum_{\ell = 3}^r\sum_{\jb \in \Ic_\ell} P_\jb(\lambda)  z_{\jb}
$$
where here the coefficients $P_\jb(\lambda)$ are unknown and where $\chi_{[\ell]}(\lambda) $ and $Z_{[\ell]}(\lambda)$ denote the homogeneous polynomials of degree $\ell$ in $\chi(\lambda)$ and $Z(\lambda)$.  

Using the assumptions on $P$, we can write
$$
P = A + B = \sum_{\ell = 3}^r P_{[\ell]} + B
$$
where $A \in \Tc_r^{\infty,\nu}$ and $B \in \Hc^s(B_s(\rho_0))$ for $s > s_0$ and $\rho_0$ sufficiently small. Moreover, $B$ has a zero of order $r+1$ at the origin. 

Identifying the coefficients of degree $\ell \leq r$ in  equation \eqref{eq:tg2}, we obtain
$$
\chi_{[\ell]}(\lambda) \circ \varphi_{H_0}^h - \chi_{[\ell]}(\lambda) = hP_{[\ell]} - h Z_{[\ell]}(\lambda) + h G_{[\ell]}(\lambda;\chi_*,P_*,Z_*) . 
$$
where $G$ is a real hamiltonian homogeneous of degree $\ell$ depending on the polynomials $\chi_{[k]}$, $P_{[k]}$ and $Z_{[k]}$ for $k < \ell$. In particular, its coefficients are polynomial of order $\leq \ell$ of the coefficients $\chi_{j}$,  $P_{j}$ and $Z_{j}$ for $j \in \Ic_k$, $k < \ell$.  

Writing down the coefficients, this equation is equivalent to 
$$
\forall\, \jb \in \Ic_r \quad
( e^{ih \Omega(\jb)} - 1)\chi_{\jb} = hP_{\jb} - h Z_{\jb} + h G_{\jb}
$$
and hence we see that the key is to control the small divisors $ e^{ih \Omega(\jb)} - 1$. 
\subsection{Non autonomous Hamiltonians}

Before giving the proof of Theorem \ref{Tmain}, we give easy results on the flow of non autonomous Hamiltonian.
Let $Q(\lambda) \in \mathcal{C}([0,1],\Tc_r^{M,\nu})$ for some $r \geq 3$, $M > 0$ and $\nu > 0$. We set
$$
\Norm{Q}{\Tc_r^{M,\nu}} = \max_{\lambda \in [0,1]} \SNorm{Q}{\Tc_r^{M,\nu}}. 
$$
The following results extend the properties already proved in \cite{Greb07} or \cite{Bam07} and needed in the proofs below. 

\begin{lemma}
\label{LP1}
Let $k \in \bar\N$, $M \in \N$, $\nu \in [0,\infty)$, $s \in \R$ with $s > \nu + 3/2$, and let $P(\lambda) \in \mathcal{C}([0,1],\Tc_{k+1}^{M,\nu})$ be a homogeneous polynomial of order $k+1$ depending on $\lambda \in [0,1]$. Then
\begin{itemize}
\item[(i)]
$P$ extends as a continuous polynomial on $\Pc_s(\C)$ depending continuously on $\lambda \in [0,1]$, and there exists a constant $C$ such that for all $z \in \Pc_s(\C)$ and all $\lambda \in [0,1]$, 
$$
|P(\lambda, z)| \leq C \Norm{P}{\Tc_{k+1}^{M,\nu}}\Norm{z}{s}^{k+1}. 
$$
\item[(ii)]
Assume moreover that $M > s+1$, then the Hamiltonian vector field $X_{P(\lambda)}$ extends as a bounded function from $\Pc_s(\C)$ to $\Pc_s(\C)$ depending continuously on $\lambda \in [0,1]$. Furthermore, for any $s > \nu +1$, there exists a constant $C$ such that for any $z \in \Pc_s(\C)$ and $\lambda \in [0,1]$, 
$$
\Norm{X_{P(\lambda)}(z)}{s} \leq C \Norm{P}{\Tc_{k+1}^{M,\nu}} \Norm{z}{s}^{k}. 
$$
\end{itemize}
\end{lemma}

\begin{lemma}
\label{LP2}
Let $r \geq 3$, $M > 0$ and let 
$$
Q(\lambda,z) = \sum_{\ell = 3}^r \sum_{\jb \in \Ic_\ell} Q_{\jb}(\lambda) z_{\jb}
$$
be an element of $\mathcal{C}([0,1],\Tc_{r}^{M,\nu})$. Let $\varphi_{Q(\lambda)}^\lambda$ be the flow associated with the non autonomous real hamiltonian $Q(\lambda)$. Then for $s > \nu +3/2$ there exist a constant $C_r$ depending on $r$ such that 
\begin{equation}
\label{Eball}
\rho < \mathrm{inf}\big(1/2, C_r \Norm{Q}{\Tc_{r}^{M,\nu}}^{-1}\big)\quad\Longrightarrow\quad
\forall\, \lambda \in [0,1], \quad \varphi_{Q(\lambda)}^\lambda(B_s(\rho)) \subset B_s(2\rho).
\end{equation}
Moreover, if $F(\lambda) \in \mathcal{C}([0,1], \Hc^s(B_s(2\rho))$ has a zero of order $r$ at the origin, then 
$F(\lambda) \circ \varphi_{Q(\lambda)}^\lambda$ has a zero of order $r$ at the origin in $B_s(\rho)$. 
\end{lemma}
\begin{Proof}
Let $z^\lambda =\varphi_{Q(\lambda)}^\lambda(z^0)$. Using the estimates of the previous lemma, we have 
$$
\begin{array}{rcl}
\displaystyle
\frac{\dd}{\dd \lambda} \Norm{z^\lambda}{s}^2 &=& 2 \langle z^\lambda , X_{Q(\lambda)}(z^\lambda) \rangle_s \\[1ex]
&\leq&  c_r \Norm{Q}{\Tc_{r}^{M,\nu}}\Norm{z^\lambda}{s}\Big(\Norm{z^\lambda}{s}^2 + \Norm{z^\lambda}{s}^{r-1}\Big)
\end{array}
$$
for some constant $c_r$ depending on $r$. 
Hence, as long as $\Norm{z^\lambda}{s} \leq 1$, we have
$$
\displaystyle\frac{\dd}{\dd \lambda} \Norm{z^\lambda}{s}^2 \leq 2 c_r \Norm{Q}{\Tc_{r}^{M,\nu}}\Norm{z^\lambda}{s}^3. 
$$
By using a comparison argument, we easily get that for $z^0 \in B_s(\rho)$ we have 
$$
\forall\, \lambda \in [0,1],\quad \Norm{z^\lambda}{s} \leq 2 \Norm{z^0}{s}. 
$$
This shows \eqref{Eball} and the rest follows. 
\end{Proof}

The next result is a consequence of Prop 6.3 in \cite{Greb07}. The only specificity is the control of the sum of the indices, and the evolution of the norm
\begin{proposition}
\label{Pbrack}
Let $k_1$ and $k_2$ two fixed integers. Let $P$ and $Q$ two homogeneous polynomials of degree $k_1+1$ and $k_2+1$ such that $P \in \mathcal{C}([0,1],\Tc_{k_1 + 1}^{M,\nu_1})$ and $Q \in \mathcal{C}([0,1],\Tc_{k_2 + 1}^{M,\nu_2})$ for some $\nu_1 > 0$, $\nu_2 > 0$ and $M > 0$. 

Then $\{P,Q\}$ defines a homogeneous polynomial of degree $k_1 + k_2$, and 
for all $M'$ and $\nu'$ such that 
$$
M' < M - \max(\nu_1,\nu_2) - 1 \quad \mbox{and} \quad \nu' > \nu_1 + \nu_2 +1, 
$$
we have $\{P,Q\} \in \mathcal{C}([0,1],\Tc_{k_1 + k_1}^{M', \nu'})$ and
$$
\Norm{\{P,Q\}}{\Tc_{k_1 + k_1}^{M', \nu'}} \leq C \Norm{P}{\Tc_{k_1 + 1}^{M,\nu_1}} \Norm{Q}{\Tc_{k_2 + 1}^{M,\nu_2}}
$$
for some constant $C$ depending on $M$, $\nu$, $M'$, $\nu'$, $k_1$ and $k_2$. 
\end{proposition}
\begin{Proof}
The proof is clear using Proposition 6.3 in \cite{Greb07}. 
We only need to verify the fact that the summations are always made over set of indices with zero moment $\Mc(\jb)$, which is trivial. 
\end{Proof} 

\begin{lemma}
Let $\chi(\lambda)$ 
be an element of $\mathcal{C}([0,1],\Tc_{r}^{M,\nu})$ for some $M >0$ and $\nu>0$. Let $\tau(\lambda) := \varphi_{\chi(\lambda)}^\lambda$ be the flow associated with the non autonomous real hamiltonian $\chi(\lambda)$. 
Let $g \in \mathcal{C}([0,1],\Tc_r^{M,\nu})$, then we can write for all $\sigma_0 \in [0,1]$, 
\begin{multline}
\label{Ereccomm}
g(\sigma_0) \circ \tau(\sigma_0) = g(\sigma_0) \\[1ex]
+\sum_{k = 0}^{r-1} \int_{0}^{\sigma_0} \cdots \int_{0}^{\sigma_k} \Big( \mathrm{Ad}_{\chi(\sigma_k)}\circ \cdots \circ  \mathrm{Ad}_{\chi(\sigma_1)} g(\sigma_0)\Big) \dd \sigma_1 \cdots \dd \sigma_k  + R(\sigma_0) 
\end{multline}
where by definition $\mathrm{Ad}_P(Q) = \{Q,P\}$
\begin{equation}
\label{ERlambda}
R(\sigma_0) =  \int_{0}^{\sigma_0} \cdots \int_{0}^{\sigma_{r}} \Big( \mathrm{Ad}_{\chi(\sigma_{r})}\circ \cdots \circ  \mathrm{Ad}_{\chi(\sigma_1)} g(\sigma_0)\Big)\circ { \tau(\sigma_{r})}\,  \dd \sigma_1 \cdots \dd \sigma_{r}. 
\end{equation}
Each term in the sum in \eqref{Ereccomm} belongs (at least) to the space $\mathcal{C}([0,1], \Tc_{kr}^{M',\nu'})$ where 
$$
\nu' = (r+1) (\nu + 2)  \quad \mbox{and}\quad  M' = M - \nu'. 
$$
The term $R(\sigma_0)$ defines an element of $\Hc^s(B_s(\rho))$ for $s > \nu'+3/2$ and $\rho \leq \mathrm{inf}(1/2,C_r \Norm{\chi}{\Tc_r^{M,\nu}}^{-1})$ and has a zero of order at least $r+1$ at the origin. 
\end{lemma}
\begin{Proof}
For a fixed $\sigma_0 \in [0,1]$, we have 
$$
\frac{\dd}{\dd \lambda} g(\sigma_0) \circ \tau(\lambda) = \{ g(\sigma_0), \chi(\lambda)\} \circ \tau(\lambda). 
$$
Hence, we have that 
$$
g(\sigma_0) \circ \tau(\sigma_0) = g(\sigma_0) + \int_{0}^{\sigma_0} \Big(\mathrm{Ad}_{\chi(\sigma_1)} g(\sigma_0)\Big) \circ \tau(\sigma_1) \, \dd \sigma_1. 
$$
Repeating again the same argument, we have
\begin{multline*}
g(\sigma_0) \circ \tau(\sigma_0) = g(\sigma_0) + \int_{0}^{\sigma_0} (\mathrm{Ad}_{\chi(\sigma_1)} g(\sigma_0)) \dd \sigma_1+\\[1ex]
 \int_{0}^{\sigma_0} \int_{0}^{\sigma_1} \Big(\mathrm{Ad}_{\chi(\sigma_2)} \circ \mathrm{Ad}_{\chi(\sigma_1)} g(\sigma_0)\Big) \circ \tau(\sigma_2) \, \dd \sigma_1 \dd \sigma_2. 
\end{multline*}
The equation \eqref{Ereccomm} is then easily shown by induction. 
The result then follows from the previous propositions. 
\end{Proof}

For a given polynomial $\chi\in  \mathcal{C}([0,1], \Tc_r^{M,\nu})$ with $r \geq 3$, we use the following notation 
\begin{equation}
\label{Echiell}
\chi(\lambda,z) = \sum_{\ell = 3}^r \chi_{[\ell]}(\lambda) = \sum_{\ell = 3}^r \sum_{\jb \in \Ic_\ell} \chi_{\jb}(\lambda) z_{\jb}
\end{equation}
where $\chi_{[\ell]}(\lambda) \in \mathcal{C}([0,1], \Tc_r^{M,\nu})$ is a homogeneous polynomial of degree $\ell$.

\begin{proposition}
\label{Pcomposition}
Let $\chi(\lambda)$ 
be an element of $\mathcal{C}([0,1],\Tc_{r}^{M,\nu})$ for some $M >0$ and $\nu>0$. Let $\varphi_{\chi(\lambda)}^\lambda$ be the flow associated with the non autonomous real hamiltonian $\chi(\lambda)$. 
Let $g \in \mathcal{C}([0,1],\Tc_r^{M,\nu})$, then we can write for all $\lambda \in [0,1]$, 
$$
g(\lambda) \circ \varphi^\lambda_{\chi(\lambda)} = S^{(r)}(\lambda) + T^{(r)}(\lambda)
$$
where 
\begin{itemize}
\item $S^{(r)}(\lambda) \in \mathcal{C}^\infty([0,1],\Tc_r^{M_1,\nu_1})$ with $\nu_1 = (r+1)(\nu +2) $ and $M_1 = M - \nu_1 $. Moreover, if we write 
$$
S(z) = \sum_{\ell = 3}^r S_{[\ell]}(\lambda)
$$
where $S_{[\ell]}(\lambda)$ is a homogeneous polynomial of degree $\ell$, then we have  for all $\ell = 3,\ldots,r$, 
$$
S_{[\ell]}(\lambda) = g_{[\ell]}(\lambda) + G_{[\ell]}(\lambda;\chi_*,g_*)
$$
where $G_{[\ell]}(\chi_*,g_*)$ is a homogeneous polynomial depending on $\lambda$ and 
the coefficients $S_j$ are polynomials of order $ < \ell$ of the coefficients appearing in the decomposition of $g$ and $\chi$. Moreover, we have
\begin{equation}
\label{Ebound1}
\Norm{G_{[\ell]}(\chi_*,g_*)}{\Tc_{r}^{M_1,\nu_1}} \leq C \Big( 1 + \sum_{m = 3}^{\ell - 1}  \Norm{g_{[m]}}{\Tc_r^{M,\nu}}^\ell\Big)\Big(1 + \sum_{m = 3}^{\ell - 1} \Norm{\chi_{[m]}}{\Tc_r^{M,\nu}}^\ell\Big)
\end{equation}
where $C$ depends on $\ell$, $M$ and $\nu$. 
\item $T^{(r)}(\lambda)\in \Hc^s(B_s(\rho))$ for $s > \nu'+3/2$ and $\rho \leq \mathrm{inf}(1/2,C_r \Norm{\chi}{\Tc_r^{M,\nu}}^{-1})$ and has a zero of order at least $r +1$ in the origin. 
Moreover, we have for all $z \in B_s(\rho)$, 
$$
\forall\, \lambda \in [0,1],\quad
\Norm{X_{T^{(r)}(\lambda)}(z)}{s} \leq C_r(\chi_*,g_*) \Norm{z}{s}^r
$$
where 
$$
C_r(\chi_*,g_*) \leq C \Big( 1 + \sum_{m = 3}^{r}  \Norm{g_{[m]}}{\Tc_r^{M,\nu}}^r\Big)\Big(1 + \sum_{m = 3}^{r} \Norm{\chi_{[m]}}{\Tc_r^{M,\nu}}^r\Big)
$$
with $C$ depending on $r$, $M$ and $\nu$. 
\end{itemize} 
\end{proposition}
\begin{Proof}
Using the previous lemma, we define $S^{(r)}$ as the polynomial part of degree less than $r$ in the expression \eqref{Ereccomm}. The remainder terms, together with the term $R(\lambda)$ in \eqref{ERlambda}, define the term $T^{(r)}(\lambda)$. The properties of $S^{(r)}(\lambda)$ and $T^{(r)}(\lambda)$ are then easily shown. 
 \end{Proof}

\subsection{Proof of the normal form result}\label{5.6}

\begin{proposition}
Let $P \in \Tc$ and $N$ be a fixed integer. 
Let $M > \nu' := (r +1 ) (\nu +2)$ and $M' = M - \nu'$. Then 
there exist 
\begin{itemize}
\item a polynomial $\chi \in \mathcal{C}([0,1], \Tc_r^{M',\nu'})$ 
$$
\chi(\lambda) = \sum_{\ell = 3}^r \chi_{[\ell]}(\lambda)  :=  \sum_{\ell = 3}^r\sum_{\jb \in \Jc_\ell(N)} \chi_\jb(\lambda)  z_{\jb}
$$

\item 
a polynomial $\chi \in \mathcal{C}([0,1], \Tc_r^{M',\nu'})$
$$
Z(\lambda) = \sum_{\ell = 3}^r Z_{[\ell]}(\lambda) := \sum_{\ell = 3}^r \sum_{\{\jb \in \Ac_\ell \cup \Ic_\ell\backslash\Jc_l(N)\}} Z_j(\lambda) z_{\jb}
$$
under normal form,\
item 
a function $R(\lambda)  \in \mathcal{C}([0,1],\Hc_s(B_s(\rho)))$ with $\rho < c_0 N^{-\beta}$ for some constant $c_0>0$ and $\beta>1$ depending on $r$, $M$, $P$, and having a zero of order $r+1$ at the origin 
\end{itemize}
 such that the following equation holds: 
\begin{equation}
\label{eq:tg22}
\forall\, \lambda \in [0,1]\quad 
  \chi(\lambda) \circ \varphi_{H_0}^h - \chi(\lambda) \circ \varphi_{hP}^{-\lambda}= hP - (h Z(\lambda) + R(\lambda)) \circ \varphi^{-\lambda}_{\chi(\lambda)} . 
\end{equation}
Furthermore there exists $C_0>0$ depending on $P$, $\nu$, $r$, $M$ such that 
$$
\left|{\chi}\right|_{\Tc_r^{M',\nu'}} + |Z|_{\Tc_r^{M',\nu'}}  \leq C_0 N^{\beta}
$$
and moreover 
$$
\forall\, \lambda \in [0,1],\quad \Norm{X_{R(\lambda)}(z)}{s} \leq C_0 N^{\beta} \Norm{z}{s}^r
$$
for $z \in B_s(\rho)$ with $\rho < c_0 N^{-\beta}$. 
\end{proposition}
\begin{Proof}
Identifying the coefficients of degree $\ell \leq r$ in the equation \eqref{eq:tg22}, we get
$$
\chi_{[\ell]} \circ \varphi_{H_0}^h - \chi_{[\ell]} = hP_{[\ell]} - h Z_{[\ell]} + h G_{[\ell]}(\chi_*,P_*,Z_*) . 
$$
where $G$ is a real hamiltonian homogeneous of degree $\ell$ depending on the polynomials $\chi_{[k]}$, $P_{[k]}$ and $Z_{[k]}$ for $k < \ell$. In particular, its coefficients are polynomial of order $\leq \ell$ of the coefficients $\chi_{\jb}$,  $P_{\jb}$ and $Z_{\jb}$ for $\jb \in \Ic_k$, $k < \ell$ and satisfying estimates of the form \eqref{Ebound1}.  
Writing down the coefficients, this equation is equivalent to 
$$
\forall\, \jb \in \Ic_r \quad
( e^{ih \Omega(\jb)} - 1)\chi_{\jb} = hP_{\jb} - h Z_{\jb} + h G_{\jb}. 
$$
We solve this equation by setting
$$
Z_\jb = P_\jb + G_\jb\quad \mbox{and}\quad \chi_\jb = 0 \quad \mbox{for}\quad \jb \in \Ac_\ell \cup \Ic_\ell\backslash\Jc_\ell(N)
$$
and
$$
Z_\jb = 0 \quad \mbox{and} \quad \chi_{\jb} = \frac{h}{e^{ih \Omega(\jb)} - 1} (P_\jb +  G_\jb) \quad \mbox{for}\quad \jb \in \Jc_\ell(N)\backslash \Ac_\ell. 
$$
Using \eqref{nonres2} and the result of Proposition \ref{Pcomposition} we  get the claimed bound for some $\beta$ depending on $r$. 

To define $R$, we simply define it by the equation \eqref{eq:tg2}. By construction, it will satisfies the announced properties. 
\end{Proof}

\begin{Proofof}{Theorem \ref{TT}}
Integrating the equation \eqref{eq:tg2} in $\lambda$, it is clear that the following equation holds: 
$$
\forall\, \lambda \in [0,1]\quad \varphi_{H_0}^h \circ \varphi_{hP}^\lambda \circ \varphi_{\chi(\lambda)}^\lambda =  \varphi_{\chi(\lambda)}^\lambda \circ\varphi_{H_0}^h \circ \varphi_{hZ(\lambda) + R(\lambda)}^\lambda.
$$
Note that using Proposition \ref{Pcomposition} and \eqref{Eball} we  show that  for $s > \nu' +1$ and $z \in B_s(\rho)$ with  $\rho = cN^{-\beta}$  we have
$$
\Norm{\varphi^\lambda_{\chi(\lambda)}(z) - z}{s} \leq C N^{\beta} \Norm{z}{s}^2.
$$
This implies in particular that 
$$
\Norm{z}{s} \leq  \Norm{\varphi^\lambda_{\chi(\lambda)}(z)}{s} + CN^{-\beta} \Norm{z}{s}
$$
For $N$ sufficiently large, this shows that $\varphi_{\chi(\lambda)}^\lambda$ is invertible and send $B_{s}(\rho)$ to $B_s(2\rho)$. 
Moreover, we have the estimate, for all $\lambda \in [0,1]$, 
$$
\Norm{\big(\varphi^\lambda_{\chi(\lambda)}\big)^{-1}(z) - z}{s} \leq C N^{\beta} \Norm{z}{s}^2.
$$
We then define $\tau = \varphi^1_{\chi(\lambda)}$ and $\psi = \varphi^1_{hZ(\lambda) + R(\lambda)}$ and verify that these application satisfy the condition of the theorem. 
\end{Proofof}

\section{Appendix: Proof of Proposition \ref{Pcrux}}

Let $\jb \in \Ic_r\backslash\Jc_r(N)$. It is clear that for $a \in \Nc$, we have, with the notation $I_a = \xi_a \eta_a$, 
$$
\{I_a,z_{\jb}\} = 0
$$
unless $(a,1)$ or $(a,-1)$ appears in $\jb$. Moreover, if this is the case, we have
$$
|\{I_a,z_{\jb}\} | \leq  2 |z_{\jb}|. 
$$
where we set 
$$
|z_{\jb}| = |z_{j_1}|\cdots |z_{j_r}|
$$
for $\jb = (j_1,\ldots,j_r) \in \Zc^r$. 
Hence we can write
$$
|\{\Bc_s^{N},Z\}(z)| \leq 2 \sum_{|k| \leq N}  |k|^{2s} \sum_{\{\jb \in \Ic_r\backslash\Jc_r(N) | \jb \supset k\} } |Z_{\jb}| |z_{\jb}|
$$
where $k = (a,\pm 1) \in \Zc$ in the first sum. 
We thus get using \eqref{Ereg}
$$
|\{\Bc_s^{N},Z\}(z)| \leq 2 \SNorm{Z}{\Tc_r^{M,\nu}} \sum_{|k| \leq N} \sum_{\{\jb \in \Ic_r\backslash \Jc_r(N) |\jb \supset k\} } |k|^{2s} \frac{\mu(\jb)^{M+\nu} }{ S(\jb)^M} |z_{\jb}|.
$$
Using Lemma \ref{L1}, the indices in the previous sum are such that at least two of them are greater than $N$. As $|k| \leq N$, these indices cannot be equal to $k$. Hence we can rewrite each indice $\jb$ containing $k$ as $(\tilde\jb,k)$ where $\tilde\jb \in \Zc^{r-1}$ contains at least two indices greater than $N$. Using the symmetries in the sum, we can moreover assume that the indices are ordered in such a way that $|j_1| > |j_2| > \ldots$. Hence, we can rewrite the previous sum as
\begin{equation}
\label{ENs1}
|\{\Bc_s^{N},Z\}(z)| \leq C \SNorm{Z}{\Tc_r^{M,\nu}} \sum_{\jb \in \Zc_{r-1},
\, |j_1|,|j_2| > N,\, |k| \leq N } |k|^{2s} \frac{\mu(\jb,k)^{M+\nu} }{ S(\jb,k)^M} |z_{\jb}| |z_k|
\end{equation}
where $\mu(\jb,k)$ and $S(\jb,k)$ denote the values of $\mu$ and $S$ associated with the $r$-uple $(j_1,j_2,\ldots,j_{r-1},k)$. Here, $C$ denote a constant depending on $r$. \\
Since $\mu(\jb,k)\leq S(\jb,k)$ and $\mu(\jb,k)\leq |j_2|$ we have for $M\geq 2$
\begin{multline*}
|\{\Bc_s^{N},Z\}(z)| \leq 
C \SNorm{Z}{\Tc_r^{M,\nu}}  \\\times \sum_{\jb \in \Zc_{r-1},\, |j_1|,|j_2| > N,\, |k| \leq N } |k|^{2s}  \left(\frac{1 }{ 1 + |j_1| - | j_2| } \right)^{2}|j_2|^{\nu +2} |z_{\jb}| |z_k|.
\end{multline*}
Then use $|k|\leq |j_1| $ to obtain
\begin{multline}
\label{Esum1}
|\{\Bc_s^{N},Z\}(z)| \leq C \SNorm{Z}{\Tc_r^{M,\nu}} \\ \times
\sum_{\jb \in \Zc_{r-1},\, |j_1|,|j_2| > N,\, |k| \leq N }  \left(\frac{1 }{ 1 + |j_1| - | j_2| } \right)^{2}|j_2|^{\nu +2} |j_1|^{s} |z_{\jb}| |k|^{s} |z_k|.
\end{multline}
By Cauchy Schwarz, one has for $s>1/2$
\begin{equation}
\label{CS}
\sum_{l\in \Zc}|z_l|\leq ||z||_s\left(\sum_{l\in \Zc} |l|^{-2s}\right)^{1/2}\end{equation}
and thus we get from \eqref{Esum1}
\begin{multline*}
|\{\Bc_s^{N},Z\}(z)| \leq C \SNorm{Z}{\Tc_r^{M,\nu}}||z||^{r-3}  \\\times \sum_{ |j_1|,|j_2| > N,\, |k| \leq N }  \left(\frac{1 }{ 1 + |j_1| - | j_2| } \right)^{2}|j_2|^{\nu +2} |j_1|^{s} |z_{\jb}| |k|^{s} |z_k|.
\end{multline*}
Hence, introducing the sequence $(b_j = |j|^{s} |z_j|)_{j \in \Zc} \in \ell^2(\Zc)$ we can write 
\begin{multline}
\label{Esum2}
|\{\Bc_s^{N},Z\}(z)| \leq C \SNorm{Z}{\Tc_r^{M,\nu}}||z||^{r-3}\\ \times 
 \sum_{ |j_1|,|j_2| > N,\, |k| \leq N }  \left(\frac{1 }{ 1 + |j_1| - | j_2| } \right)^{2}|j_2|^{\nu +2-s}b_{j_2} b_{j_1}b_k.
\end{multline}
Moreover, the sum in $|k| \leq N$ in  \eqref{Esum2} yields by Cauchy-Schwartz inequality
$$
\sum_{|k| \leq N} b_k \leq C N^{d/2} \sqrt{\Bc_s^{N}(z)} \leq C N^{d/2} \Norm{z}{s}. 
$$
where  $d$ is the dimension of $\Nc = \N^d$ or $\Z^d$. 
Hence, we get from  \eqref{Esum2} using $ |j_2| > N$
$$
|\{\Bc_s^{N},Z\}(z)| \leq C  N^{-s+2+\nu + d/2} \SNorm{Z}{\Tc_r^{M,\nu}}||z||^{r-2} \sum_{ |j_1|\geq|j_2| > N}  \left(\frac{1 }{ 1 + |j_1| - | j_2| } \right)^{2}b_{j_2} b_{j_1}
$$
and this concludes the proof of \eqref{ENsP} since, if $a$ and $c$ are two sequences in $\ell^2(\Zc)$ we have by a convolution argument
\begin{equation}
\label{convol}
\sum_{ j,l}  \left(\frac{1 }{ 1 + |j| - | l| } \right)^{2}|a_{l}|\ | c_{j}|\leq C  ||a||_{\ell^2(\Zc)}\ ||c||_{\ell^2(\Zc)}
\end{equation}
for some universal constant $C$.

 Note that \eqref{ENsP2} is an easily shown by similar calculations (the only difference lies in the fact that there is no summation in $|k| \leq N$). 

We now show \eqref{ERsP}.\\
As $\Rc^N_s$ contains only indices greater than $N$, we can write (see \eqref{ENs1})
\begin{equation}\label{estimR}
|\{\Rc_s^{N},Z\}(z)| \leq C \SNorm{Z}{\Tc_r^{M,\nu}} \sum_{\jb \in \Ic_{r-1},\, |j_1| > N,\, |k| > N } |k|^{2s} \frac{\mu(\jb,k)^{M+\nu} }{ S(\jb,k)^M} |z_{\jb}| |z_k|
\end{equation}
where the sum is made over ordered indices $|j_1| > |j_2|> \cdots$. Note that in opposition with the previous situation, we cannot ensure that $|j_2| > N$ in this sum.  
We first notice  that, for all $k$ and $\jb$,
\begin{equation}
\label{Emagic}
|k| \frac{\mu(\jb,k) }{S(\jb,k)} \leq 2 |j_1|,
\end{equation}
Actually,  if $|k|\leq 2j_{1}$ then \eqref{Emagic} holds true since $
\frac{\mu(\jb,k)}{S(\jb,k)}\leq 1$. Now if $k\geq 2j_{1}$
then $S(l,j)\geq ||k|
-|j_{1}||\geq \frac{1}{2} |k|$ and thus
$$
|k|\frac{\mu(\jb,k)}{S(\jb,k)}\leq 2{\mu(\jb,k)}\leq
2|j_{1}| .$$
Then we distinguish two cases in this sum \eqref{estimR}: 
$$
|\{\Rc_s^{N},Z\}(z)| \leq C \SNorm{Z}{\Tc_r^{M,\nu}}(I_1 + I_2 )
$$
corresponding to the two cases $|j_2 | \leq |k| $, ($I_1$) and  $|k| < |j_2|$, ($I_2$). 

\medskip
\textbf{Case 1: $|j_2| \leq |k|$} 

In this situation, we use \eqref{Emagic}, $\mu(\jb,k)= |j_2|$, $\mu(\jb,k)\leq S(\jb,k)$ to conclude for $M\geq s+2$
$$
I_1\leq 2^s \sum_{\jb \in \Ic_{r-1},\, |j_1| > N,\, |k| > N } |k|^{s}|j_1|^s\left(\frac{ 1}{1 + ||j_1| - |k||} \right)^{2}|j_2|^{2+\nu} |z_{\jb}| |z_k|.
$$
Then use \eqref{CS} and the notation $(b_j = |j|^{s} |z_j|)_{j \in \Zc} \in \ell^2(\Zc)$ to get
\begin{align*}
I_1&\leq 2^s||z||_s^{r-3} \sum_{j_2, |j_1| > N,\, |k| > N } \left(\frac{ 1}{1 + ||j_1| - |k||} \right)^{2}|j_2|^{2+\nu-s}b_{j_2} b_k b_{j_1}\\
&\leq C ||z||_s^{r-2}\Rc_s^N(z)
\end{align*}
where we have used again \eqref{convol} for $\sum_{j_1,k}$ and \eqref{CS} for $\sum_{j_2}$ .

\medskip
\textbf{Case 2: $|j_2| \geq |k| $} 

In this situation, we still have $\mu(\jb,k) \leq |j_2|$ and using that both $|j_1|$ and $|j_2|$ are greater than $|k|$ we get for $M\geq 2$
\begin{align*}
I_2 &\leq C \Norm{z}{s}^{r-3} \sum_{|j_1| > N, \, |j_2 | \geq |k| > N} |j_1|^s  |j_2|^{s/2} |k|^{s/2}
\left(\frac{1}{1 + |j_1| - |j_2|}\right)^{2}|j_2|^{2+\nu}  |z_{j_1} |z_k| |z_{j_2}| \\
&\leq 
 C \Norm{z}{s}^{r-3} \sum_{|j_1| > N, \, |j_2 | \geq |k| > N}
\left(\frac{1}{1+ |j_1| - |j_2|}\right)^{2} b_{j_1}\frac{b_k}{|k|^{s/2}} \frac{b_{j_2}}{|j_2|^{s/2-2-\nu}} \\
&\leq 
 C \Norm{z}{s}^{r-3} \Rc_s^N(z)^{3/2}
\end{align*}
where in the last inequality, we used that $ b_{j_1}$, $\frac{b_k}{|k|^{s/2}}$ and $\frac{b_{j_2}}{|j_2|^{s/2-2-\nu}} $ are respectively in $\ell^2(\Zc)$, $\ell^1(\Zc)$ (for $s>1$) and $\ell^2(\Zc)$ (for $s\geq 4+2\nu$) and we used again \eqref{convol} and \eqref{CS}.


\begin{thebibliography}{30}

\bibitem{Bam03}
D.~Bambusi, \emph{Birkhoff normal form for some nonlinear {PDE}s}, Comm. Math.
  Physics \textbf{234} (2003), 253--283.
  
\bibitem{Bam07}
D.~Bambusi,
 \emph{A birkhoff normal form theorem for some semilinear pdes},
  Hamiltonian Dynamical Systems and Applications, Springer, 2007, pp.~213--247.
  
\bibitem{BDGS}
D.~Bambusi, J.-M. Delort, B.~Gr{\'e}bert, and J.~Szeftel, \emph{Almost global
  existence for {H}amiltonian semilinear {K}lein-{G}ordon equations with small
  {C}auchy data on {Z}oll manifolds}, Comm. Pure Appl. Math. \textbf{60}
  (2007), no.~11, 1665--1690. 
  
\bibitem{BG06}
{\rm D. Bambusi and B. Gr\'ebert},
{\em Birkhoff normal form for PDE's with tame modulus}. Duke Math. J.  135  no. 3 (2006), 507–-567.


\bibitem{CHL08a}
{\rm D. Cohen, E. Hairer and C. Lubich}, 
{\em Long-time analysis of nonlinearly perturbed wave equations via modulated Fourier expansions},
Arch. Ration. Mech. Anal. 187 (2008) 341-368. 

\bibitem{DF07}
{\rm G. Dujardin and E. Faou},
{\em Normal form and long time analysis of splitting schemes for the linear Schrödinger equation with small potential.}
Numerische Mathematik 106, 2 (2007) 223--262 

\bibitem{FGP1}
{\rm E. Faou, B. Gr\'ebert and E. Paturel},
{\em Long time analysis of splitting methods applied to discretized Hamiltonian partial differential equations.}

\bibitem{CHL08b}
{\rm  E. Hairer and C. Lubich}, 
{\em Spectral semi-discretisations of weakly nonlinear wave equations over long times}, 
Found. Comput. Math. 8 (2008) 319-334.

\bibitem{CHL08c}
{\rm D. Cohen, E. Hairer and C. Lubich},
{\em Conservation of energy, momentum and actions in numerical discretizations of nonlinear wave equations}, 
Numerische Mathematik 110 (2008) 113--143.

\bibitem{GL08a}
{\rm L. Gauckler and C. Lubich},
{\em Nonlinear Schr\"odinger equations and their spectral discretizations over long times},
Preprint (2008).

\bibitem{GL08b}
{\rm L. Gauckler and C. Lubich},
{\em Splitting integrators for nonlinear Schr\"odinger equations over long times},
Preprint (2008). 


\bibitem{Greb07}
{\rm B. Gr\'ebert},  
{\em Birkhoff normal form and Hamiltonian PDEs.}
S\'eminaires et Congr\`es 15 (2007), 1--46
 
 \bibitem{GIP}
{\rm B. Gr\'ebert, E. Paturel and R. Imekraz},  
{\em Long time behavior for solutions of   semilinear Schr\"odinger equation with harmonic potential and small Cauchy data on $\R^d$.}
Preprint (2008)

\bibitem{DS06} 
{\rm J. M. Delort and J. Szeftel}, 
{\em Long-time existence for semi-linear Klein-Gordon equations with small cauchy data on Zoll manifolds}, Amer. J. Math
128 (2006), 1187--1218.
 
\bibitem{Deuf79}
{\rm P. Deuflhard}
{\em A study of extrapolation methods based on multistep schemes without parasitic solutions}.
Z. angew. Math. Phys. 30 (1979) 177-189. 

\bibitem{HLW}
{\rm E. Hairer, C. Lubich and G. Wanner}
{\em Geometric Numerical Integration. Structure-Preserving Algorithms for Ordinary Differential Equations}. Second Edition. Springer 2006. 


\bibitem{Shan00}
{\rm Z. Shang}
{\em Resonant and Diophantine step sizes in computing invariant tori of Hamiltonian systems}
Nonlinearity 13 (2000), 299--308.


\end{thebibliography}
\end{document}